\documentclass[a4paper,12pt, reqno]{amsart}
\usepackage{amssymb,amsmath,latexsym}
\usepackage{color}%cyan,red;magenta,green;yellow,blue
\usepackage{comment}
\usepackage{hyperref}
\allowdisplaybreaks

\newcommand{\normal}{\color{black}}
\newcommand{\bk}{\color{black}}

\def\1{{\bf 1}}

\def\nn{\nonumber}
\def\bee{\begin{equation}}
\def\eee{\end{equation}}
 \def\sB {{\mathcal B}}

\def\R {{\mathbb R}}

\newtheorem{thm}{Theorem}[section]
\newtheorem{lemma}[thm]{Lemma}
\newtheorem{defn}[thm]{Definition}
\newtheorem{prop}[thm]{Proposition}
\newtheorem{corollary}[thm]{Corollary}

\numberwithin{equation}{section}

\def\qed{{\hfill $\Box$ \bigskip}}
\def\NN{{\mathcal N}}

\def\FF{{\mathcal F}}
\def\EE{{\mathcal E}}

\def\R{{\mathbb R}}

\def\E{{\mathbb E}}

\def\P{{\mathbb P}}
\def\N{{\mathbb N}}
\def\eps{\varepsilon}
\def\wh{\widehat}
\def\wt{\widetilde}
\def\pf{\noindent{\bf Proof.} }

\usepackage{hyperref}
\begin{document}
\title[
Potential theory of Dirichlet forms degenerate at the boundary]
{Potential theory of
Dirichlet forms degenerate at the boundary: the case of no killing potential}
\author{ Panki Kim \quad Renming Song \quad and \quad Zoran Vondra\v{c}ek}
\thanks{Panki Kim: This work was  supported by the National Research Foundation of
Korea(NRF) grant funded by the Korea government(MSIP) (No. NRF-2015R1A4A1041675)
}
\thanks{Renming Song: Research supported in part by a grant from
the Simons Foundation (\#429343, Renming Song)}
\thanks{Zoran Vondra\v{c}ek: Research supported in part by the Croatian Science Foundation under the project 4197.}

 \date{}

\begin{abstract}
In this paper we consider the Dirichlet form on the half-space $\R^d_+$ defined by the jump kernel $J(x,y)=|x-y|^{-d-\alpha}\sB(x,y)$, where $\sB(x,y)$ can be degenerate at the boundary. 
Unlike our previous works \cite{KSV, KSV21} where we imposed critical killing,
here we assume that the killing potential is identically zero. 
In case $\alpha\in (1,2)$ we first show that the corresponding Hunt process has finite lifetime and dies at the boundary. Then, as our main contribution, we prove 
 the boundary Harnack principle and establish sharp two-sided Green function estimates. 
Our results cover the case of the censored $\alpha$-stable process,  $\alpha\in (1,2)$, in the half-space studied in \cite{BBC}.
\end{abstract}
\maketitle

\bigskip
\noindent {\bf AMS 2020 Mathematics Subject Classification}: Primary 60J45; Secondary 60J50, 60J76.

\bigskip\noindent
{\bf Keywords and phrases}: Jump processes, 
jump kernel, 
jump kernel degenerate at the boundary,
Carleson estimate, boundary Harnack principle, Green function.

\smallskip
%%%%%%%%%%%%%%%%%%%%%%%%%%%%%%%%%%%%%%%%%%%%%%%%%%%%%%%%%%%%%%%%%%%%%%%%%%%%%%%%%%%%%%
%%%%%%%%%%%%%%%%%%%%                                            Introduction                                                                               %%%%%%%%%%%%%%%%%%%%%
%%%%%%%%%%%%%%%%%%%%%%%%%%%%%%%%%%%%%%%%%%%%%%%%%%%%%%%%%%%%%%%%%%%%%%%%%%%%%%%%%%%%%%

\section{Introduction}

Let $\R_+^d=\{x=(\wt{x},x_d):\, x_d>0\}$ be the upper half-space in the $d$-dimensional Euclidean space $\R^d$. In this paper we study the Dirichlet form $(\EE, \FF)$ on $L^2(\R^d_+, dx)$ defined by
\begin{equation}\label{e:form}
\EE(u,v):=\frac12 \int_{\R^d_+}\int_{\R^d_+} (u(x)-u(y))(v(x)-v(y))J(x,y)\, dy\, dx,
\end{equation}
where $\FF$ is the closure of $C_c^{\infty}({\R^d_+})$ under $\EE_1:=\EE+(\cdot, \cdot)_{L^2({\R^d_+},dx)}$. Our main assumption is 
on
the jump kernel $J(x,y)$: We assume that $J(x,y)=|x-y|^{-d-\alpha}\sB(x,y)$, $\alpha\in (0,2)$, where $(x,y)\mapsto \sB(x,y)$ is a symmetric function satisfying  certain H\"older-type and scaling conditions, and most importantly, is comparable to the function
\begin{eqnarray}\label{e:B(x,y)}
\wt{B}(x,y)&:=& \Big(\frac{x_d\wedge y_d}{|x-y|}\wedge 1\Big)^{\beta_1}\Big(\frac{x_d\vee y_d}{|x-y|}\wedge 1\Big)^{\beta_2} \left[ \log\Big(1+\frac{(x_d\vee y_d)\wedge |x-y|}{x_d\wedge y_d\wedge |x-y|}\Big)\right]^{\beta_3}\nn \\  
 & &   \times  \left[\log \Big(1+\frac{|x-y|}{(x_d\vee y_d)\wedge |x-y|}\Big)\right]^{\beta_4}. \normal
\end{eqnarray}
Here $\beta_1, \beta_2, \beta_3, \beta_4$ are non-negative parameters such that $\beta_1>0$ if $\beta_3>0$, and $\beta_2>0$ if $\beta_4>0$. 
Here and below,  $a\wedge b:=\min \{a, b\}$, $a\vee b:=\max\{a, b\}$.
The precise assumptions on $\sB(x,y)$ are given in Section 2. Although we allow that $\sB(x,y)\equiv 1$, our focus is on the case when $\beta_1 \vee \beta_2>0$. In such a case, the function $\sB(x,y)$ vanishes at the boundary of $\R^d_+$, and we call the corresponding Dirichlet form degenerate at the boundary. We refer to  $\sB(x,y)$ as the boundary part of the jump kernel $J(x,y)$. This setting was introduced in \cite[Section 5]{KSV}. The Hunt process  associated with the Dirichlet form $(\EE, \FF)$ will be denoted by $Y=(Y_t, \P_x)$ and its lifetime by $\zeta$.

Our motivation to study the form \eqref{e:form} and the corresponding process $Y$ comes from two sources.

Firstly, note that in case $J(x,y)=|x-y|^{-d-\alpha}$ (i.e.~ $\sB(x,y)\equiv 1$), the process $Y$ is the censored $\alpha$-stable process in the half-space $\R^d_+$ which 
was
introduced and studied in \cite{BBC} (also for a more general state space than $\R^d_+$). Two main results of \cite{BBC} can be roughly described as follows: (1) There is a dichotomy between cases $\alpha\in (1,2)$ and $\alpha\in (0,1]$. In the former case the process $Y$ has finite lifetime $\zeta$ and approaches the boundary of the state space at $\zeta$, while in the latter, $Y$ is conservative and will never approach the boundary; (2)
In case when the state space $D$ is a $C^{1,1}$ open set and $\alpha\in (1,2)$,  the boundary Harnack principle holds with the exact decay rate $\delta_D(x)^{\alpha-1}$ (here $\delta_D(x)$ denotes the distance of the point $x$ to the boundary of $D$). Shortly after, in case of a bounded $C^{1,1}$ open set and $\alpha\in (1,2)$,  sharp two-sided Green function estimates were established in \cite{CK02}. 

Secondly, a Dirichlet form related to \eqref{e:form} was introduced in \cite{KSV} and further studied in \cite{KSV21}. We considered the form $(\EE^{\kappa}, \FF^{\kappa})$ where
\begin{equation}\label{e:form-killing}
\EE^{\kappa}(u,v)=\EE(u,v)+\int_{\R^d_+}u(x)v(x)\kappa(x)\, dx\, ,
\end{equation}
and $\FF^{\kappa}=\FF\cap L^2(\R^d_+, \kappa(x)dx)$. The killing function is given by $\kappa(x)=C(\alpha, p, \sB)x_d^{-\alpha}$, where $C(\alpha, p, \sB)$ is a semi-explicit strictly positive and finite constant depending on $\alpha$, $\sB$ and a parameter $p\in ((\alpha-1)_+, \alpha+\beta_1)$. The investigation of the form \eqref{e:form-killing} was 
initiated  in \cite{KSV} 
and completed in \cite{KSV21} with two main results: Sharp two-sided Green function estimates for all admissible values of the parameters involved in $\wt{B}(x,y)$, cf.~\cite[Theorem 1.1]{KSV21}, and full identification of the parameters for which the boundary Harnack principle holds true, cf.~\cite[Theorem 1.2 and Theorem 1.3]{KSV21}. In proving those results, the strict positivity of the killing function was used in an essential way in several places. This includes the proof of finite lifetime, Carleson estimate, and the decay of the Green function at the boundary. 

The goal of this paper is to extend the main results of \cite{KSV21} to the Dirichlet form \eqref{e:form} (which has no killing) in case $\alpha\in (1,2)$. Due to the fact that  $\lim_{p\downarrow (\alpha-1)_+}C(\alpha, p, \sB)=0$, this can be considered as a limiting case of the setting in \cite{KSV21}. Theorem \ref{t:BHP} below can be viewed as a generalization of the corresponding result in \cite{BBC} in case of the state space $\R^d_+$ and $\alpha\in (1,2)$, to jump kernels degenerate at the boundary, while  Theorem \ref{t:Green} is related to the main result of \cite{CK02}.

The following two theorems are the main contribution of this paper. 
For $a,b>0$ let $D_{\wt{w}}(a,b):=\{x=(\wt{x}, x_d)\in \R^d:\, |\wt{x}-\wt{w}|<a, 0<x_d<b\}$.
Assumptions \textbf{(A1)}-\textbf{(A4)} are given in Section \ref{s:SP}.
\begin{thm}\label{t:BHP}
Suppose  that $\alpha\in (1,2)$
and  that $\sB$ satisfies \normal \textbf{(A1)}-\textbf{(A4)}. There exists $C \ge 1$ such that for all $r>0$, $\wt{w} \in \R^{d-1}$,  and any non-negative function $f$ in $\R^d_+$ which is harmonic in $D_{\wt{w}}(2r, 2r)$ with respect to $Y$ and vanishes continuously on $B(({\wt{w}},0), 2r)\cap \partial \R^d_+$,  we have 
\begin{equation}\label{e:TAMSe1.8new}
\frac{f(x)}{x_d^{ \alpha-1}}\le C\frac{f(y)}{y_d^{\alpha-1}}, \quad x, y\in D_{\wt{w}}(r/2, r/2) .
\end{equation}
\end{thm}

We would like to mention here that, even though the rate of the boundary decay of harmonic functions is different than the one in \cite{KSV}, we still use the main result of \cite{KSV} to prove  Theorem \ref{t:BHP}. More specifically,  to control the exit distributions  in Lemma \ref{e:POTAe7.14} from below, 
we use the  lower bound for the corresponding exit distributions in \cite{KSV}. 
The latter is the most  complicated and technical part in  \cite{KSV}.
We point out in passing the following observation: 
to establish some boundary theory for non-local operators with singular kernels (with no critical killing),  the corresponding 
non-local operators with critical killing can play an important role.

Let $G(x,y)$, $x,y\in \R^d_+$, denote the Green function of the process $Y$ (see Section \ref{s:EGP} for the existence of the Green function).
\begin{thm}\label{t:Green}
Suppose that $\alpha\in (1,2)$ and $d > (\alpha+\beta_1+\beta_2)\wedge 2$.
Assume that $\sB$ satisfies \textbf{(A1)}-\textbf{(A4)}.
Then there exists $C>1$ such that for all $x,y\in \R^d_+$,
\begin{eqnarray}\label{e:Green}
\lefteqn{C^{-1} \left(\frac{x_d}{|x-y|}  \wedge 1 \right)^{\alpha-1}\left(\frac{y_d}{|x-y|}  \wedge 1 \right)^{\alpha-1} \frac{1}{|x-y|^{d-\alpha}} \le G (x,y)} \nonumber \\
& \le & C \left(\frac{x_d}{|x-y|}  \wedge 1 \right)^{\alpha-1}\left(\frac{y_d}{|x-y|}  \wedge 1 \right)^{\alpha-1} \frac{1}{|x-y|^{d-\alpha}}\, .
\end{eqnarray}
\end{thm}

Note that in both results we have assumed that $\alpha\in (1,2)$. The case $\alpha\in (0,1]$  is qualitatively different and new methods are needed to analyze it. 
We leave this case for future research.

Now we explain the content of the paper, our strategy of proving the results and differences to the methods used in \cite{BBC} and \cite{KSV, KSV21}.

In Section \ref{s:SP} we precisely introduce the setup and assumptions on the boundary function $\sB(x,y)$, and recall some of the relevant results from \cite{KSV}. 

The goal of Section \ref{s:H} is to prove that in case $\alpha\in (1,2)$, the process $Y$ has finite lifetime and is therefore transient. 
The proof is new and relies on a
Hardy-type inequality, see Proposition \ref{t:hardy}. This inequality implies that $\FF\neq \overline{\FF}$, where $\overline{\FF}$ is the closure of $C_c^{\infty}({\overline \R^d_+})$ under $\EE_1=\EE+(\cdot, \cdot)_{L^2({\R^d_+},dx)}$. 
This implies that $Y$ is a
(proper) subprocess of $\overline{Y}$ -- the Hunt process associated with $(\EE, \overline{\FF})$, hence the lifetime of $Y$ is finite. A consequence of finite lifetime is  Corollary \ref{c:lemma4-1} which has two parts: The first one shows that the process $Y$ approaches the boundary at the lifetime, while the second part
\normal replaces \cite[Lemma 4.1]{KSV} in the standard proof of the Carleson inequality, see Theorem \ref{t:carleson}. 

Section \ref{s:D} is devoted to proving Dynkin's formula for some non-compactly supported and non-smooth functions. Let
\begin{equation}\label{e:operator}
L_{\alpha}^\sB f(x):=\textrm{p.v.}\int_{\R^d_+}(f(y)-f(x))J(x,y)\, dy\
\end{equation}
be the operator corresponding to the form $(\EE, \FF)$, defined for all $f:\R^d_+\to \R$ for which the principal value integral makes sense. It is straightforward to see that $L_{\alpha}^\sB x_d^{\alpha-1}=0$, which can be understood as $x\mapsto x_d^{\alpha-1}$ being harmonic in the analytic sense. See Lemma \ref{l:LB-on-g} below. In order to use probabilistic methods, it is crucial to show that this function is harmonic in the probabilistic sense. 
The proofs of \cite[Lemmas 3.3 and  5.1]{BBC} rely
on using the isotropic stable process and 
its part process in $\R^d_+$. 
Since these two processes are of no help to us in the present setting,
we use instead  Dynkin's formula for barriers, cf.~Proposition \ref{l:dynkin-hp}. The proof of this formula is a slight modification of 
the arguments in \cite[Section 9]{KSV}. 

Section \ref{s:BHP} is devoted to the proof of Theorem \ref{t:BHP}. We first argue that 
the proofs of some results from \cite{KSV} in the case $p>\alpha-1$ are easily modified to the case $p=\alpha-1$. 
Then we show that for any function $f$ as in the statement of Theorem \ref{t:BHP}, it holds that
$$
\frac{f(x)}{f(y)}\asymp \frac{ \P_x\big(Y_{\tau_{D_{\wt{w}}(r/2, r/2)}}\in D(1, 1)\big)} { \P_y\big(Y_{\tau_{D_{\wt{w}}(r/2, r/2)}}\in D(1, 1)\big)}, \quad x,y\in D_{\wt{w}}(r/2, r/2).
$$
Since $x\mapsto x_d^{\alpha-1}$ satisfies the conditions in Theorem \ref{t:BHP},
the assertion of Theorem \ref{t:BHP} is valid.

In the first part of Section  \ref{s:EGP} we present the proof of Theorem \ref{t:Green}. The proof uses some results from \cite{KSV21}, scaling and the boundary Harnack principle. In the second part we give sharp estimates of the Green potential of 
$x_d^\gamma$ for $\gamma>-\alpha$.
Again, we argue that, using the boundary Harnack principle, proofs of some lower bounds of the killed Green function  obtained in \cite{KSV21} for $p>\alpha-1$ are valid without any change for the case $p=\alpha-1$. Having these Green function estimates one can apply results from \cite{AGV} to get the estimates of the Green potentials.

We end the introduction with an explanation of  the connection between the process $Y$ and the process $Y^{\kappa}$ associated with the Dirichlet form $(\EE^{\kappa}, \FF^{\kappa})$. 
This connection is analogous to the one between the censored stable process and the killed stable process, cf.~\cite[Theorem 2.1]{BBC}. Namely, the process $Y$ can be obtained from $Y^{\kappa}$ through either the Ikeda-Nagasawa-Watanabe piecing together procedure, or through the Feynman-Kac transform via $\exp\int_0^t \kappa(Y^{\kappa}_t)dt$. The case $\sB\equiv 1$ and $\kappa(x)=C(\alpha/2, \alpha, \sB)x_d^{-\alpha}$ corresponds exactly to the isotropic $\alpha$-stable process killed upon exiting $\R^d_+$.

Throughout this paper, the positive constants $\beta_1$, $\beta_2$, $\beta_3$, $\beta_4$,
$\theta$ , $r_0$, $n_0$ will remain the same.
We will use the following convention:
Lower case letters 
$c, c_i, i=1,2,  \dots$ are used to denote constants in the proofs
and the labeling of these constants starts anew in each proof.
The notation $c_i=c_i(a,b,c,\ldots)$, $i=0,1,2,  \dots$ indicates  constants depending on $a, b, c, \ldots$. We will not specify the dependency on $d$.
We will use ``$:=$" to denote a
definition, which is read as ``is defined to be".
For any $x\in \R^d$ and $r>0$, we use $B(x, r)$ to denote the open ball of radius $r$ centered at $x$.

\bigskip

\section{Setup and Preliminary}\label{s:SP}
In this section we precisely describe the setup and recall some preliminary results from earlier works. 

Let $d \ge 1$, $\alpha\in (0,2)$, $j(|x-y|)=|x-y|^{-\alpha-d}$ and $J(x,y)=j(|x-y|)\sB(x,y)$. We first give the assumptions on the  boundary function $\sB(x,y)$. \normal

\noindent
\textbf{(A1)} $\sB(x,y)=\sB(y,x)$ for all $x,y\in  \R^d_+$.

\medskip
\noindent
\textbf{(A2)}  
 If $\alpha \ge1$, \bk there exist $\theta>\alpha-1$ and  
$C>0$ such that 
$$
|\sB(x, x)-\sB(x,y)|\le 
C\left(\frac{|x-y|}{x_d\wedge y_d}\right)^{\theta}\,.
$$ 

\medskip
\noindent
\textbf{(A3)}
There exist $C\ge 1$ and 
parameters $\beta_1, \beta_2, \beta_3,   \beta_4 \normal  \ge 0$,  
with $\beta_1>0$ if $\beta_3 >0$,  and $\beta_2>0$ if $\beta_4>0$, 
such that
\begin{equation}\label{e:B7}
C^{-1}\wt{B}(x,y)\le \sB(x,y)\le C \wt{B}(x,y)\, ,\qquad x,y\in \R^d_+\, ,
\end{equation}
where
$\wt{B}(x,y)$ is defined in \eqref{e:B(x,y)}.

\medskip
\noindent
\textbf{(A4)} 
For all $x,y\in \R^d_+$ and  $a>0$, $\sB(ax,ay)=\sB(x,y)$. 
In case $d\ge 2$, for  all 
$x,y\in \R^d_+$ and $\wt{z}\in \R^{d-1}$, $\sB(x+(\wt{z},0), y+(\wt{z},0))=\sB(x, y)$.

For examples of functions $\sB$ satisfying  \textbf{(A1)}-\textbf{(A4)}, see \cite{KSV, KSV21}.
 Assumption \normal  \textbf{(A3)} implies that $\sB(x, y)$ is bounded. 
Note that  \textbf{(A4)} implies that $x\mapsto \sB(x, x)$ is a constant on $\R^d_+$. Without loss
of generality, we will assume that $\sB(x, x)=1$.

 We observe that if $\beta_4>0$, 
  then, for any $\eps\in (0, \beta_2)$,  there exists $c_\eps>0$ such that  
  \begin{align}\label{e:B7_2} 
 &(\log 2)^{-\beta_4}\wt{B}_{\beta_1, \beta_2, \beta_3, 0}(x,y) \le \wt{B}_{\beta_1, \beta_2, \beta_3, \beta_4}(x,y) \le c_{\eps} \wt{B}_{\beta_1, \beta_2-\eps, \beta_3, 0}(x,y).
 \end{align}

{\it 
 Throughout the paper 
we always assume that}
\begin{align*}
  J(x,y)&=j(|x-y|)
\sB(x,y) \text{ on }  \R_+^d\times \R_+^d \text{ with  } \sB
\text{ satisfying } \textbf{(A1)}-\textbf{(A4)} \text{ and } \sB(x, x)=1.
\end{align*}

 Let $\overline \R_+^d=\{x=(\wt{x},x_d):\, x_d \ge 0\}$.
Define
\begin{align*}
\EE(u,v)&:=\frac12 \int_{\R^d_+}\int_{\R^d_+} (u(x)-u(y))(v(x)-v(y))J(x,y)\, dy\, dx\\
&=
\frac12 \int_{\overline \R^d_+}\int_{\overline \R^d_+} (u(x)-u(y))(v(x)-v(y))J(x,y)\, dy\, dx.
\end{align*}
By Fatou's lemma,  $(\EE, C_c^{\infty}({\R^d_+}))$ and $(\EE, C_c^{\infty}({\overline \R^d_+}))$
are closable in 
$L^2({\R^d_+}, dx)(=L^2({\overline \R^d_+}, dx))$.
Let $\FF$ be the closure of $C_c^{\infty}({\R^d_+})$ under
$\EE_1:=\EE+(\cdot, \cdot)_{L^2({\R^d_+},dx)}$ and 
let $\overline \FF$ be the closure of $C_c^{\infty}({\overline \R^d_+})$ under
$\EE_1=\EE+(\cdot, \cdot)_{L^2( \R^d_+, dx)}$. 
Then $(\EE, \FF)$ and $(\EE, \overline \FF)$
are regular Dirichlet 
forms.

Let $((Y_t)_{t\ge 0}, (\P_x)_{x\in {\R^d_+}\setminus \NN})$ be the Hunt 
process associated with  $(\EE, \FF)$ whose  lifetime is $\zeta$. 
By \cite[Proposition 3.2]{KSV},  the exceptional set 
$\NN$ can be taken to be the empty set. 
We add a cemetery point $\partial$ to the state space ${\R^d_+}$ and define $Y_t=\partial$
for $t\ge \zeta$. 
Let $((\overline Y_t)_{t\ge 0}, (\P_x)_{x\in {\overline \R^d_+}\setminus \NN_0})$ be the Hunt 
process associated with  $(\EE, \overline \FF)$ where $\NN_0$ is  an exceptional set.
Let $(\EE, \overline{\FF}_{\R^d_+})$ be the part form of 
 $(\EE, \overline{\FF})$ on $\R^d_+$.
i.e.,  the form corresponding to the process $\overline{Y}$ killed at 
the exit time
$\tau_{\R^d_+}:=\inf\{t>0:\, \overline{Y}_t\notin \R^d_+\}$. It follows from \cite[Theorem 4.4.3(i)]{FOT} that 
$(\EE, \overline{\FF}_{\R^d_+})$
is a regular Dirichlet form on $L^2({\R^d_+},dx)$ and that $C^{\infty}_c(\R^d_+)$ is its core. Hence $\overline{\FF}_{\R^d_+}=\FF$ implying that $\overline{Y}$ killed upon exiting $\R^d_+$ is equal to $Y$. Thus   we conclude that $Y$ is a subprocess of $\overline{Y}$, 
that the exceptional set 
$\NN_0$ can be taken to be a subset of $\partial \R^d_+,$ 
and
that the lifetime of $Y$ can be identified with $\tau_{\R^d_+}$. 

Suppose that for all $x\in \R^d_+$ it holds that $\P_x(\tau_{\R^d_+}=\infty)=1$. Then $(Y_t, \P_x, x\in \R^d_+)\stackrel{d}{=} (\overline{Y}_t, \P_x, x\in \R^d_+)$ implying that $\FF=\overline{\FF}_{\R^d_+}=\overline{\FF}$.

 For any $r>0$, define a process $Y^{(r)}$ by
$Y^{(r)}_t:=r Y_{r^{-\alpha} t}$.
By the proof of \cite[Lemma 5.1]{KSV}, $Y$ has the following scaling property. 
\begin{lemma}\label{l:scaling-of-Y}
 $(Y^{(r)},
  \P_{x/r})\normal $ has the same law as $(Y, \P_x)$.
\end{lemma}

For any open subset $V$ of $\R^d_+$ and for $r>0$, we define 
$rV:=\{rx:\, x\in V\}$ and $\tau_V=\inf\{t>0:\, Y_t\notin V\}$.
A consequence of Lemma \ref{l:scaling-of-Y} is that 
\begin{equation}\label{e:exit-time-scaling}
\E_{rx}\tau_{rV}=r^{\alpha}\E_x \tau_V\, , \qquad x\in V\, .
\end{equation}

\begin{defn}\label{D:1.1}  \rm 
A non-negative Borel function defined on
$\R_+^d$ is said to be  {harmonic}
in an open set $V\subset \R_+^d$ with respect to
$Y$ if for every bounded open set $U\subset\overline{U}\subset V$,
\begin{equation}\label{e:har}
f(x)= \E_x \left[ f(Y_{\tau_{U}})\right] \qquad
\hbox{for all } x\in U.
\end{equation}
A non-negative Borel function $f$ defined on $\R_+^d$ is said to be \emph{regular harmonic} in an open set $V\subset \R_+^d$ if 
$$
f(x)= \E_x \left[ f(Y_{\tau_{V}})\right] \qquad \hbox{for all } x\in V.
$$
\end{defn}

The following result is taken form \cite{KSV, KSV21}.

\begin{thm}[Harnack inequality, {\cite[Theorem 1.1]{KSV}} \& {\cite[Theorem 1.4]{KSV21}}] \label{t:uhp}
\begin{itemize}
\item[(a)] 
There exists a constant 
$C_1>0$
such that for any
$r>0$, any $B(x_0, r) \subset \R^d_+$ and any non-negative function $f$ in $\R^d_+$ which is 
harmonic in $B(x_0, r)$ with respect to $Y$, we have
$$
f(x)\le C_1 
f(y), \qquad \text{ for all } x, y\in  B(x_0, r/2).
$$

\item[(b)] 
There exists a constant 
$C_2>0$ 
such that for any  $L>0$, any  $r>0$,
any $x_1,x_2 \in \R^d_+$ with $|x_1-x_2|<Lr$ and $B(x_1,r)\cup B(x_2,r) \subset \R^d_+$ 
and any non-negative function $f$ in $\R^d_+$ which is  harmonic in $B(x_1,r)\cup B(x_2,r)$ with respect to $Y$, we have
$$
f(x_2)\le C_2
(L+1)^{\beta_1+\beta_2+d+\alpha} f(x_1)\, .
$$
\end{itemize}
\end{thm}
\section{Hardy inequality and the finite lifetime}\label{s:H}

For a given  $\sB$, we define $C(\alpha, p, \sB)$ for $\alpha \in (0,2)$ and $p\in (-1, \alpha+\beta_1)$ by
\begin{equation}\label{e:explicit-C}
C(\alpha, p, \sB)=
\int_{\R^{d-1}}\frac{1}{(|\wt{u}|^2+1)^{(d+\alpha)/2}}
\int_0^1 \frac{(s^p-1)(1-s^{\alpha-p-1})}{(1-s)^{1+\alpha}}
\sB\big((1-s)\wt{u}, 1), s\mathbf{e}_d \big)\, ds
d\wt{u}\, ,
\end{equation}
where $\mathbf{e}_d=(\tilde{0}, 1)$.
In case $d=1$, $  C(\alpha, p, \sB)$ is defined by
$$
C(\alpha, p, \sB)=\int_0^1 \frac{(s^p-1)(1-s^{\alpha-p-1})}{(1-s)^{1+\alpha}}
\sB\big(1, s \big)\, ds, 
$$
but we will only give the statement of the result for $d\ge 2$. The statement in the $d=1$ case is similar and simpler.

We first note that 
$p\mapsto C(\alpha, p, \sB)$ is 
strictly increasing for $p\in ((\alpha-1)_+, \alpha+\beta_1)$ (see \cite[Lemma  5.4  and Remark  5.5]{KSV}) and $\lim_{p\uparrow \alpha+\beta_1}C(\alpha, p, \sB)=\infty$.
Moreover,
\begin{align}
\label{e:C}
C(\alpha, p, \sB) 
\begin{cases}
\in (0, \infty) & \text{ for } p\in ((\alpha-1)_+,  \alpha+\beta_1);\\
=0 & \text{ for } p=0, \alpha-1;\\
\in (-\infty, 0) & \text{ for } p\in (\alpha-1, 0) \cup (0, \alpha-1).
\end{cases}
\end{align}

Let 
$$
C_c^2(\R^d_+; \R^d)=\{f:\R^d_+\to  \R : \text{ there exists } u \in  C_c^2(\R^d) \text{ such that } u=f \text{ on } \R^d_+\}
$$
be the space of functions on $\R^d_+$ that are restrictions of $C_c^2(\R^d)$ functions. Clearly, if $f\in C_c^2(\R^d_+; \R^d)$ then $f\in C^2_b(\R^d_+) \cap L^2(\R^d_+)$. 

For $\eps>0$, let 
\begin{align}\label{e:defn-LBDe}
L_{\alpha, \eps}^\sB f(x):=\int_{\R^d_+, |y-x|> \varepsilon}(f(y)-f(x))J(x,y)\, dy, \quad x\in \R^d_+,
\end{align}
so that 
\begin{align}\label{e:defn-LBD}
L_{\alpha}^\sB f(x)=\textrm{p.v.}\int_{\R^d_+}(f(y)-f(x))J(x,y)\, dy=
\lim_{\varepsilon \to 0} L_{\alpha, \eps}^\sB f(x),
\end{align}
which is defined for all functions $f:\R^d_+\to \R$ for which the principal value integral makes sense. 
We have shown in \cite[Proposition 3.4]{KSV} that this is the case when $f\in C_c^2(\R^d_+; \R^d)$.

 For $p\in \R$,  let $g_{p}(x)=x_d^p$,  $x \in \R^d_+$. 

\begin{lemma}\label{l:LB-on-g}
For  $p\in (-1, \alpha+\beta_1)$, it holds that
\begin{align}\label{e:Lp}
L_{\alpha}^\sB g_{p}(x)= C(\alpha,p, \sB)x_d^{p-\alpha}, \quad x\in \R^d_+.
\end{align}
In particular, 
\begin{align}\label{e:Lal}
L_{\alpha}^\sB g_{\alpha-1}(x)=0, \quad x\in \R^d_+.
\end{align}
Moreover, there exists $\wh C=\wh C(\alpha,p, \sB)>0$ such that 
\begin{align}\label{e:Lpe}
|L_{\alpha, \eps}^\sB g_{p}(x)| \le  \wh C \,x_d^{p-\alpha} \quad \text{for all }x\in \R^d_+ \text{ and } \eps \in  (0, x_d/2].
\end{align}
\end{lemma}
\pf  The equality 
\eqref{e:Lp} is proved in \cite[Lemma 5.4]{KSV}  for $p\in ((\alpha-1)_+, \alpha+\beta_1)$. 
 It is easy to see that the proof in fact works for $p\in (-1, \alpha+\beta_1)$. 
We now follow the proof in \cite[Lemma 5.4]{KSV} to show \eqref{e:Lpe}.

Fix $x=(\wt{0}, x_d) \in \R^d_+$ and $\eps \in  (0, x_d/2]$.
By the change of variables $y=x_d z$,  and by using \textbf{(A4)}, 
we  have
\begin{align*}
L_{\alpha}^\sB g_{p, \eps}(x)
&=x_d^{p-\alpha}  \int_{\R^d_+, |\wt{z}|^2+| z_d-1|^2> (\varepsilon/x_d)^2}\frac{z_d^p-1}{| (\wt{z}, z_d)-\mathbf{e}_d|^{d+\alpha}}\sB(\mathbf{e}_d,  (\wt{z}, z_d))\, 
dz_d d \wt{z}\\
&
=:x_d^{p-\alpha} I_1(\eps)\, .
\end{align*}
Using the change of variables $\wt{z}=|z_d-1|\wt{u}$, we get
\begin{align*}
I_1(\eps)&=
 \int_{\R^d_+, | z_d-1|^2 |\wt{u}|^2+| z_d-1|^2> (\varepsilon/x_d)^2}
{(|\wt{u}|^2+1)^{-(d+\alpha)/2}}\frac{z_d^p-1}{|z_d-1|^{1+\alpha}} 
\sB\big(\mathbf{e}_d, (|z_d-1|\wt{u}, z_d)\big)dz_d \, d\wt{u}\\
&=
 \int_{\R^{d-1}}{(|\wt{u}|^2+1)^{-(d+\alpha)/2}}
 I_2(\eps, \wt{u})
 \, d\wt{u},
\end{align*}
where
 $$I_2(\eps, \wt{u})=
 \left(\int_0^{1-(\eps/x_d)(|\wt{u}|^2+1)^{-1/2} }+\int_{1+(\eps/x_d)(|\wt{u}|^2+1)^{-1/2}}^{\infty}\right) 
\frac{z_d^p-1}{|z_d-1|^{1+\alpha}} 
\sB\big(\mathbf{e}_d, (|z_d-1|\wt{u}, z_d)\big)dz_d .
$$
Fix $\wt{u}$ and let $\epsilon_0=(\eps/x_d)(|\wt{u}|^2+1)^{-1/2} \le 1/2$. By
the same argument as that in  the proof of \cite[Lemma 5.4]{KSV}, we have that
$
I_2(\eps, \wt{u})=I_{21}(\eps, \wt{u})+I_{22}(\eps, \wt{u})
$
where
\begin{align*}
I_{21}(\eps, \wt{u})&:=\int_0^{1-\epsilon_0}\frac{(s^p-1)+(s^{\alpha-1-p}-s^{\alpha-1})}{(1-s)^{1+\alpha}}\sB\big(((1-s)\wt{u}, 1), s\mathbf{e}_d \big)\, ds, \\
I_{22}(\eps, \wt{u}) &:=\int_{1-\epsilon_0}^{\frac{1}{1+\epsilon_0}}\frac{s^{\alpha-1-p}-s^{\alpha-1}}{(1-s)^{1+\alpha}}\sB\big(((1-s)\wt{u}, 1), s\mathbf{e}_d \big)\, ds,
\end{align*}
and there exists a constant $c_2>0$ independent of $\wt{u}\in \R^{d-1}$ such that
\begin{align}
\label{e:I21n}
|I_{21}(\eps, \wt{u})|  \le \int_0^1 \frac{|(s^p-1)(1-s^{\alpha-p-1})|}{(1-s)^{1+\alpha}}\sB\big(((1-s)\wt{u}, 1), s\mathbf{e}_d \big)\, ds<c_1 <\infty\,  .
\end{align}
Moreover, by \cite[p.121]{BBC} and the fact that $\epsilon_0 \le 1/2$,
$$
\left|\int_{1-\epsilon_0}^{\frac{1}{1+\epsilon_0}}\frac{s^{\alpha-1-p}-s^{\alpha-1}}{(1-s)^{1+\alpha}} \sB\big((1-s)\wt{u}, 1), s\mathbf{e}_d \big)ds\right| \le c_2 \epsilon_0^{2-\alpha}  \le c_2.
$$
Therefore,
\begin{align*}
 \sup_{\eps \in  (0, x_d/2]}|I_1(\eps)|& \le c_3\int_{\R^{d-1}}\frac{d\wt{u}}{(|\wt{u}|^2+1)^{(d+\alpha)/2}} =c_4<\infty\, .
\end{align*}
\qed

 We now show that the following Hardy inequality holds when $\alpha \not= 1$. \normal
\begin{prop}\label{t:hardy} 
Suppose $\alpha \not= 1$.
Then there exists $C=C(\alpha) \in (0, \infty)$ 
 such that 
for all  $u \in \FF,$
\begin{align}
\label{e:hardy}
\EE(u,u) \ge C \int_{\R^d_+} \frac{u(x)^2}{x_d^\alpha}dx.
\end{align}

\end{prop}

\pf
Since $\FF$ is the closure of $C_c^{\infty}({\R^d_+})$ under
$\EE_1$, it suffices to prove \eqref{e:hardy} for $u \in C_c^{\infty}({\R^d_+})$.

Fix $u \in C_c^{\infty}({\R^d_+})$, 
choose a $p\in (\alpha-1, 0) \cup (0, \alpha-1)$ and let $v(x)=u(x)/g_p(x)$.  
Recall from \eqref{e:C} that $C(\alpha, p, \sB) \in  (-\infty, 0)$.

Using the elementary identity 
$
(ab-cd)^2=a^2b(b-d)+c^2d(d-b)+bd(a-c)^2$ and the symmetry of $J$, we have that, for all $ \eps >0$,
\begin{align*}
& \int_{\R^d_+ \times \R^d_+, |x-y|>\eps} (v(y)g_p (y)-v(x)g_p (x))^2J(x,y)\, dy\, dx\\
=&  \int_{\R^d_+ \times \R^d_+, |x-y|>\eps} 
 v(y)^2g_p (y) (g_p (y)-g_p (x)) + v(x)^2g_p (x) (g_p (x)-g_p (y)) 
 J(x,y)\, dy\, dx\\&+
 \int_{\R^d_+ \times \R^d_+, |x-y|>\eps} 
  g_p (x)g_p (y)(v(y)-v(x))^2
J(x,y)\, dy\, dx\\
 \ge& -2 \int_{\R^d_+ }  v(x)^2g_p (x)\left(\int_{ \R^d_+, |x-y|>\eps} 
 (g_p (y)-g_p (x)) J(x,y)\, dy\right)\, dx\\
 =& -2\int_{\text{supp}(u) }  \frac{u(x)^2}{g_p (x)} L_{\alpha, \eps}^\sB g_{p}(x)
  dx.
 \end{align*}
Let $a_0:=$dist$(\R^d_{-},$supp$(u))/2>0$. By \eqref{e:Lpe},
the functions $\{ \frac{u(x)^2}{g_p (x)} L_{\alpha, \eps}^\sB g_{p}(x): \eps\in (0, a_0) \}$ are uniformly bounded on 
\text{supp}$(u)$. Thus, by the bounded convergence theorem,  \eqref{e:C} and \eqref{e:Lp}, 
\begin{align*}
&\EE(u,u) 
=\lim_{\eps \downarrow 0} \frac{1}2\int_{\R^d_+ \times \R^d_+, |x-y|>\eps} (v(y)g_p (y)-v(x)g_p (x))^2J(x,y)\, dy\, dx\\
&\ge -\lim_{\eps \downarrow 0} \int_{\text{supp}(u) }  \frac{u(x)^2}{g_p (x)} L_{\alpha, \eps}^\sB g_{p}(x)
dx=  c \int_{\R^d_+ }  \frac{u(x)^2}{g_p (x)}x_d^{p-\alpha} dx=c \int_{\R^d_+ }  \frac{u(x)^2}{x_d^\alpha}dx,
\end{align*}
where $c=- C(\alpha, p, \sB) \in (0, \infty)$.
\qed

Recall that  $\zeta$ is the lifetime of $Y$.
Using the above 
Hardy inequality, we now show that  $\zeta$ is finite when $ \alpha>1$.

\begin{prop}
\label{p:finitelife}
Suppose $ \alpha>1$. Then
$\FF \not=\overline \FF$ and $\P_x(\zeta<\infty)=1$ for all $x \in \R^d_+$.
\end{prop}
\pf 
Take a $u\in  C_{ c }^{\infty}(\overline \R^d_+)$ 
such that $u \ge 1$ on $B(0, 1) \cap \R^d_+$, then 
$u \notin \FF$. In fact, if $u \in \FF$, then by Proposition \ref{t:hardy}, 
$$
\infty > \EE(u,u) \ge c  
\int_{\R^d_+ }  \frac{u(x)^2}{x_d^\alpha}dx \ge c 
\int_{B(0, 1) \cap \R^d_+ }  |x|^{-\alpha}dx=\infty, 
$$
which gives a contradiction. 

The fact that  $\FF \not=\overline \FF$ implies that there is a point $x_0  \in \R^d_+$ such that $\P_{x_0}(\zeta<\infty)>0$. 
Then by the scaling property of $Y$ in Lemma \ref{l:scaling-of-Y}, we have that $\P_x(\zeta<\infty)=\P_{x_0}(\zeta<\infty)>0$ for all $x \in \R^d_+$.
Now, by the same argument as in the proof of \cite[Proposition 4.2]{BBC}, we have that  $\P_x(\zeta<\infty)=1$ for all $x \in \R^d_+$.
\qed

The fact that the lifetime of $Y$ is finite has two important consequences.

\begin{corollary}\label{c:lemma4-1}
\begin{itemize}
	\item[(a)] For all $x\in \R^d_+$, $\P_x(Y_{\zeta-}\in \partial \R^d_+)=1$.
	\item[(b)] There exists a constant $n_0 \ge 2$ such that for all $x\in \R^d_+$, $\P_x\left(  \tau_{B(x,n_0 x_d)}=\zeta     \right)> 1/2$. 
\end{itemize}
\end{corollary}
\pf Using Lemma \ref{l:scaling-of-Y}, we see that 
$$
\P_x\big(  \tau_{B(x,n x_d)}=\zeta\big)=  \P_{(\wt 0, 1)}\big(  \tau_{B((\wt 0, 1) ,n )} =\zeta  \big), \quad x \in \R^d_+.
$$
The sequence of events $(\{ \tau_{B((\wt 0, 1) ,n )}=\zeta \})_{n\ge 1}$ 
is increasing in $n$ and 
\begin{align}\label{e:zetaf}
\cup_{n=1}^\infty\big\{  \tau_{B((\wt 0, 1) ,n )}=\zeta  \big \}=\big\{\zeta<\infty\big\}.
\end{align}
Thus, by Proposition \ref{p:finitelife} we have
\begin{equation}\label{e:c-lemma4-1}
\lim_{n \to \infty} \P_{(\wt 0, 1)}\big(  \tau_{B((\wt 0, 1) ,n )}=\zeta  \big)=
\P_{(\wt 0, 1)}\big(\zeta<\infty\big)=1.
\end{equation} 
Moreover, since there is no killing inside $\R^d_+$, it holds that $\{ \tau_{B((\wt 0, 1) ,n )}=\zeta\}\subset \{Y_{\zeta-}\in \partial \R^d_+\}$ for each $n\ge 1$. Thus it follows from \eqref{e:zetaf} and \eqref{e:c-lemma4-1} that $\P_{(\wt 0, 1)}(Y_{\zeta-}\in \partial \R^d_+)=1$. The claim (a) now follows by scaling.

To see (b), note that by \eqref{e:c-lemma4-1} there exists a $n_0 \ge 2$ such that $\P_{(\wt 0, 1)}\big(  \tau_{B((\wt 0, 1) ,n_0 )} =\zeta  \big)>1/2$.
Therefore,  
$$
\P_x\big(  \tau_{B(x,n_0 x_d)}=\zeta    \big)=\P_{(\wt 0, 1)}\big(  \tau_{B((\wt 0, 1) ,n_0 )} =\zeta  
\big)>1/2, \quad x \in \R^d_+.
$$
\qed

\section{Dynkin's formula for barriers}\label{s:D}

In this section we always assume that $\alpha>1$.

 Recall that \normal for $a,b>0$ and $\wt{w} \in \R^{d-1}$,
$$
D_{\wt{w}}(a,b):=\{x=(\wt{x}, x_d)\in \R^d:\, |\wt{x}-\wt{w}|<a, 0<x_d<b\}.
$$
Without loss of generality, we will mostly deal with the case $\wt w=\wt{0}$. 
We will write $D(a,b)$ for $D_{\wt{0}}(a,b)$ and
$U(r)=D_{\wt{0}}(\frac{r}2, \frac{r}2).$
Further  we use  $U$ for $U(1)$.
For any $a>0$, set $D(a):=\{x=(\widetilde{x}, x_d)\in \R^d: x_d>a\}$ and $U_a(r):= \{y\in U(r): \delta_{U(r)}>a\}$. We write $U_a$ for $U_a(1)$.

Let $v\in C^{\infty}_c(\R^d)$ be a non-negative smooth radial function such that $v(y)=0$ for $|y|\ge 1$ and $\int_{\R^d}v(y)\, dy=1$. 
For $b\ge 10$ and $k\ge 1$,  set $v_k(y):=b^{kd} v(b^ky)$. Next we define
$g_k:=v_k\ast(g{\bf 1}_{D(5^{-k})})$ for a bounded, compactly supported function $g$ vanishing on $\R^d\setminus  \R^d_+$.
Since $b^{-k}<5^{-k}$, we have $g_k\in C^\infty_c(\R^d_+)$ and hence $L_{\alpha}^\sB g_k$ is defined everywhere. Also note that $v_k\ast g \in C^{\infty}_c(\R^d_+; \R^d)$ and thus $L_{\alpha}^\sB (v_k \ast g)$ is well defined (cf. \cite[Subsection 3.2]{KSV}). 

Let $(a_k)_{k\ge 1}$ be a decreasing sequence of positive numbers such that $\lim_{k\to \infty}a_k=0$ and  
$$
a_k\ge 2^{-k(\beta_1/2+1)/(1+\alpha+3\beta_1/2)}\ge 2^{-k}.
$$

\begin{lemma}\label{l:Lvk}
Let $R, M \ge 1$ and $g:\R^d \to [0,M]$ be a bounded, compactly supported function vanishing on $\R^d\setminus  \R^d_+$. For any $z\in U(R)$, it holds that
\begin{equation}\label{e:Lvk-1}
\lim_{k\to \infty}L_{\alpha}^\sB(v_k\ast g -g_k)(z)=0\, .
\end{equation}
Moreover,  there exists $C>0$ independent of $R, M \ge 1$ and $g$ such that  for all $k \ge 2$ and $z\in U_{a_k}(R)$,
\begin{equation}\label{e:Lvk-2}
0\le L_{\alpha}^\sB(v_k\ast g -g_k)(z) \le CM (2/3)^{k(\beta_1/2+1)} z_d^{\beta_1}\, .
\end{equation}
\end{lemma}
\pf 
Let  $z\in U_{a_k}(R)$. We first estimate the difference 
\begin{eqnarray*}
L_{\alpha}^\sB(v_k\ast g -g_k)(z)&=&\lim_{\epsilon \to 0}\int_{\R^d_+, |y-z|>\epsilon}\frac{((v_k\ast g)(y)-g_k(y))-((v_k\ast g)(z)-g_k(z))}{|y-z|^{d+\alpha}}\sB(y,z)\, dy.
\end{eqnarray*}
Note that for $k\ge 2$, $u\in B(0, b^{-k})$ and $y\in \R^d_+$ with $y_d>3^{-k}$, 
it holds that
$y_d-u_d>3^{-k}-10^{-k}>5^{-k}$. Therefore
\begin{align}\label{e:3.11}
1-{\bf 1}_{D(5^{-k})}(y-u)=0.
\end{align}
Since $v_k$ is supported in $B(0, b^{-k})$, for all $k\ge 2$ and $z\in \R^d_+$ with $z_d> a_k>2^{-k}$,
$$
\int_{\R^d}(1-{\bf 1}_{D(5^{-k})}(z-u))g(z-u)v_k(u)du=0.
$$
Thus $(v_k\ast g -g_k)(z)=0$. 
Due to the same reason we have that for $z\in U_{a_k}(R)$,
\begin{align*}
&\int_{\R^d_+, |y-z|>\epsilon}\frac{\big((v_k\ast g)(y)-g_k(y)\big)-\big((v_k\ast g)(z)-g_k(z)\big)}{|y-z|^{d+\alpha}}\sB(y,z)\, dy\\
& =\int_{\R^d_+, |y-z|>\epsilon, y_d\le 3^{-k}}\int_{\R^d} v_k(u)
\frac{(1-{\bf 1}_{D(5^{-k})})(y-u)g(y-u)}{|y-z|^{d+\alpha}}\, du\sB(y,z)\, dy \\
&\le  M  \int_{\R^d}v_k(u)\, du \int_{\R^d, y_d\le 3^{-k}}\frac{\sB(y,z)}{|y-z|^{d+\alpha}}\, dy\\
&\le  c_1M\int_{y_d\le 3^{-k}}\frac{1}{|y-z|^{d+\alpha}}\left(\frac{y_d}{|y-z|}\right)^{\beta_1/2}\, dy\\
&\le  c_2M (3^{-k})^{\beta_1/2+1}\int_0^{\infty}\frac{t^{d-2}}{(t^2+cz_d^2)^{(d+\alpha+\beta_1/2)/2}}\, dt \\
&=  c_3M (3^{-k})^{\beta_1/2+1} z_d^{-1-\alpha-\beta_1/2} \int_0^{\infty}\frac{s^{d-2}}{(s^2+1)^{(d+\alpha+\beta_1/2)/2}}\, ds \\
&\le  c_4M (2/3)^{k(\beta_1/2+1)} z_d^{\beta_1}\, .
\end{align*}
In the third line 
we used that $0\le g\le M$,  
in the fourth the fact that (together with \eqref{e:B7_2}) 
$$
\left(\frac{y_d\wedge z_d}{|y-z|}\wedge 1\right)^{\beta_1}\log\left(1+\frac{(y_d\vee z_d)\wedge|y-z|}{y_d\wedge z_d \wedge|y-z|}\right)^{\beta_3}
\le c \left(\frac{y_d}{|y-z|}\wedge 1\right)^{\beta_1/2},
$$
in the fifth integration in polar coordinates in $\R^{d-1}$, in the sixth the change of variables $t=c^{1/2}z_d s$, and in the last line the fact that $2^{-k(\beta_1/2+1)}z_d^{-1-\alpha-3\beta_1/2}\le 1$ which follows from $z_d\ge a_k$ and the choice of $a_k$.  Note also that it is clear from the second line that the first line is non-negative. Thus by letting $\epsilon \to 0$ we get for $z\in U_{a_k}(R)$,
$$
0\le L_{\alpha}^\sB(v_k\ast g -g_k)(z) \le c_4  M (2/3)^{k(\beta_1/2+1)} z_d^{\beta_1}\, .
$$
Now take $z\in U(R)$. Then there exists $k_0\ge 1$ such that $z\in U_{a_k}(R)$ for all $k\ge k_0$, and it follows from above that
$$
\lim_{k\to \infty}L_{\alpha}^\sB(v_k\ast g -g_k)(z)=0\, .
$$
\qed

\begin{lemma}\label{l:vLg}
Assume that $R, M \ge 1$ and $g:\R^d_+\to [0,M]$  is a function which is $C^2$  on $D(R,R)$. 
For any $k\ge 2$, $z\in U_{a_k}(R)$ and $|u|<b^{-k}$,
\begin{equation}\label{e:vLg-0}
\mathrm{p.v.} \int_{\R^d_+} \frac{g(y-u)-g(z-u)}{|y-z|^{d+\alpha}}\sB(y,z)\, dy
\end{equation}
is well defined. Moreover, for $z\in U_{a_k}(R)$,
\begin{equation}\label{e:vLg-1}
L_{\alpha}^\sB(v_k\ast g)(z)=\int_{\R^d}v_k(u)\left( \mathrm{p.v.}  \int_{\R^d_+} \frac{g(y-u)-g(z-u)}{|y-z|^{d+\alpha}}\sB(y,z)\, dy\right)\, du\, ,
\end{equation}
and there exists 
 $C(z)=C(z, g, M, R)>0$ such that $|L_{\alpha}^\sB (v_k\ast g)(z)|\le C(z)$ for all $k\ge 2$.
\end{lemma}
\pf Let $z\in U_{a_k}(R)$ and $|u|<b^{-k}$. 
Let $G(y, z, u):=(g(y-u)-g(z-u))|y-z|^{-d-\alpha}$.
For $0<\epsilon <\eta <z_d/10$, consider
\begin{align*}
&\int_{\R^d_+, \epsilon < |y-z| }G(y, z, u)\sB(y,z)\, dy-\int_{\R^d_+, \eta < |y-z| }G(y, z, u)\sB(y,z)\, dy\\
&=\int_{\R^d_+, \epsilon < |y-z| < \eta}G(y, z, u)\sB(y,z)\, dy\\
&= \int_{ \epsilon < |y-z| <\eta}G(y, z, u)\, dy+ \int_{ \epsilon < |y-z| <\eta}G(y, z, u)(\sB(y,z)-1)\, dy\\
&=: I+II\, .
\end{align*}
Since $g$ is $C^2$ on $D(R,R)$ and $y-u, z-u\in D(R,R)$, we see that
\begin{align*}
&| I |\le \int_{\epsilon < |(y-u)-(z-u)|<\eta} \frac{|g(y-u)-g(z-u)-\nabla g(z-u){\bf 1}_{(|(y-u)-(z-u)|<1)}\cdot (y-z)|}{|(y-u)-(z-u)|^{d+\alpha}}\, dy\\
&\le c_1 \sup_{w\in B(z, z_d/5)}|\partial^2 g (w)| \int_{ \epsilon < |y-z| <\eta} |y-z|^{-d-\alpha+2}\, dy=c_2(z) (\eta^{2-\alpha}-\epsilon^{2-\alpha})\, .
\end{align*}
Further, by using the mean value theorem in the first line and  
\textbf{(A2)}  in the second, we get 
\begin{align*}
&| II |\le  \sup_{w\in B(z, z_d/5)}|\nabla g (w)| \int_{\epsilon < |y-z| <\eta} \frac{|\sB(y,z)-1|}{|y-z|^{d+\alpha-1}}\, dy\\
&\le c_3(z)  \int_{\epsilon < |y-z| <\eta} |y-z|^{-d-\alpha}\left(\frac{|y-z|}{y_d\wedge z_d}\right)^{\theta}\, dy\\
&\le c_4 c_3(z) z_d^{-\theta} \int_{\epsilon < |y-z| <\eta} |y-z|^{-d-\alpha+\theta}\, dy=c_5(z) (\eta^{\theta-\alpha+1}-\epsilon^{\theta-\alpha+1})\, .
\end{align*}
The estimates  for $I$ and $II$ imply that the principal value integral in \eqref{e:vLg-0} is well defined.

Let $z\in U_{a_k}(R)$. For $\epsilon < z_d/10$ and $|u|<b^{-k}$, we have
\begin{align*}
&\left| \int_{\R^d_+, |y-z|>\epsilon} G(y, z, u)\sB(y,z)\, dy\right|\\
&\le \left|\int_{\R^d_+,  |y-z|\ge z_d/10} G(y, z, u)\sB(y,z)\, dy\right| +\left|\int_{\R^d_+, \epsilon < |y-z| < z_d/10} G(y, z, u)\sB(y,z)\, dy\right| \\
&=:III+IV\, .
\end{align*}
Estimating $g$ by $M$,   we get that
$$
III\le  2M  \int_{|y-z|\ge z_d/10}|y-z|^{-d-\alpha}dy \le c_7 z_d^{-\alpha}=c_8(z) \, .
$$
The integral in $IV$ is estimated in $I$ and $II$ with $\eta=z_d/10$, so we have 
$$
IV\le c_2(z)(z_d/10)^{2-\alpha}+ c_5(z)(z_d/10)^{\theta-\alpha +1}= c_9(z)\, .
$$
Thus we have that
\begin{equation}\label{e:vLg-2}
\left| \int_{\R^d_+, |y-z|>\epsilon} \frac{g(y-u)-g(z-u)}{|y-z|^{d+\alpha}}\sB(y,z)\, dy\right| \le c_{10}(z)\, .
\end{equation}
Hence we can use the dominated convergence theorem to conclude that
\begin{align*}
&L_{\alpha}^\sB (v_k\ast g)(z) =\lim_{\epsilon \to 0}\int_{\R^d_+, |y-z|>\epsilon}\frac{(v_k\ast g)(y)-(v_k\ast g)(z)}{|y-z|^{d+\alpha}}\sB(y,z)\, dy\\
&=\lim_{\epsilon\to 0} \int_{|u|<b^{-k}}v_k(u) \int_{\R^d_+, |y-z|>\epsilon} \frac{g(y-u)-g(z-u)}{|y-z|^{d+\alpha}}\sB(y,z)\, dy \, du\\
&=\int_{|u|<b^{-k}} v_k(u)\left(\lim_{\epsilon \to 0} \int_{\R^d_+, |y-z|>\epsilon} \frac{g(y-u)-g(z-u)}{|y-z|^{d+\alpha}}\sB(y,z)\, dy\right)\, du\, ,
\end{align*}
which is \eqref{e:vLg-1}. The last statement follows from \eqref{e:vLg-2} . 
\qed

We note also that 
if $g$ is continuous in $D(R, R)$, then $\lim_{k\to \infty}(v_k\ast g)(z)=g(z) $ for all $z\in D(R, R)$.
Let $h_{p, R}(x)=x_d^p {\bf 1}_{D(R,R)}(x)$  and   $h_{p, \infty}(x)=g_{p}(x)=x_d^p$ 
for $x \in \R^d_+$. We also let $h_{p}(x)=h_{p, 1}(x)$.
\begin{lemma}\label{l:hq}
  Assume that $R \ge 1$. \bk
Let $p\in [\alpha-1, \alpha+\beta_1)$ and set  $b=10\vee 2^{4(p-2)_-+3}$ and   \bk
$$
a_k:=2^{-k(p+1+\frac12\beta_1)/(\alpha+1+\frac32\beta_1)}\vee 2^{-k(2+\beta_1)/(1+\alpha+\frac32\beta_1-p)}.
$$
There exists a constant $C=C(R)>0$ such that
for any  $k\ge 1$ and $z\in U_{a_k}(R)$,
\begin{equation}\label{e:hq-1}
|L_{\alpha}^\sB(v_k \ast h_{p, R})(z)-L_{\alpha}^\sB h_{p, R}(z)|\le C \left(\frac{4}{5}\right)^k \, .
\end{equation}
In particular, the functions 
$z\mapsto |L_{\alpha}^\sB(v_k \ast h_{p, R})(z)-L_{\alpha}^\sB h_{p, R}(z)|$ 
are all bounded
by the constant $C$on $U(R)$, and for any $z\in U(R)$, 
$
\lim_{k\to \infty} \left| L_{\alpha}^\sB(v_k \ast h_{p, R})(z)-L_{\alpha}^\sB h_{p, R}(z)\right| =0.
$
\end{lemma}
\pf First note that $a_k\ge 2^{-k}$ since the first term in its definition is larger than $2^{-k}$. 
 Fix $ z \ \in U_{a_k}( R \bk)$. 
By using
Lemma \ref{l:vLg} with $g=h_{p, R}$ in the second   line below we see that
\begin{align*}
&L_{\alpha}^\sB(v_k\ast h_{p, R})(z)\\
&=\int_{\R^d}v_k(u)\left(\mathrm{p.v.}\int_{\R^d_+}\frac{h_{p, R}(y-u)- h_{p, R}(z-u)}{|y-z|^{d+\alpha}}\sB(y, z)dy\right)du\\
&=\int_{\R^d}v_k(u)\left(\mathrm{p.v.}\int_{\R^d_+}\frac{h_{p, R}(y-u)- h_{p, R}(z-u)-(h_{p, R}(y)-h_{p, R}(z))}{|y-z|^{d+\alpha}}\sB(y, z)dy\right)du\\
&\quad+ L_{\alpha}^\sB h_{p, R}(z).
\end{align*}

 Set $b=10\vee 2^{4(p-2)_-+3}$. Now we write, for $u\in B(0, b^{-k})$,
\begin{align*}
&
 \mathrm{p.v.} \bk \int_{\R^d_+}\frac{h_{p, R}(y-u)- h_{p, R}(z-u)-(h_{p, R}(y)-h_{p, R}(z))}{|y-z|^{d+\alpha}}\sB(y, z)dy\\
&=\int_{D(R+b^{-k}, R+b^{-k})\setminus U(R), y_d>5^{-k}}+\int_{D(R+b^{-k}, R+b^{-k}), y_d<5^{-k}}\\
&\quad + 
\int_{U(R), y_d>5^{-k}, |y-z|>2^{-1}z_d}+\ 
  \mathrm{p.v.} \int_{U(R) \cap B(z, z_d/2)}\bk=:I+II+III+IV.
\end{align*}

We deal with $I$ first. For $u\in B(0, b^{-k})$,
\begin{align*}
I&=\int_{D(R+b^{-k}, R+b^{-k})\setminus D(R-b^{-k}, R-b^{-k}), y_d>5^{-k}}+
\int_{D(R-b^{-k}, R-b^{-k})\setminus U(R), y_d>5^{-k}}=:I_1+I_2.
\end{align*}
Obviously, we have
$
|I_1|\le c_1(R) b^{-k}.
$

Let $A_k:=(D(R-b^{-k}, R-b^{-k})\setminus U(R))\cap\{y:y_d>5^{-k}\}$  and
$$F_1(y_d, z_d, u_d):=p(y_d-z_d)\cdot \int^1_0\left((z_d-u_d+t(y_d-z_d))^{p-1}-(z_d+t(y_d-z_d))^{p-1}\right)dt.$$ \bk
Then, we have
\begin{align*}
&|I_2|=\big|\int_{A_k}
\frac{(y_d-u_d)^p-y_d^p-((z_d-u_d)^p-z_d^p)}{|y-z|^{d+\alpha}}
\sB(y, z)dy\big|\\
& =\big|\int_{A_k}\frac{
 F_1(y_d, z_d, u_d) 
}{|y-z|^{d+\alpha}}
\sB(y, z)dy\big|
\le c|u_d| \int_{A_k}\frac{(z_d \wedge y_d)^{-(p-2)_-} |y_d-z_d|}{|y-z|^{d+\alpha}}dy\\
&
\le c b^{-k}5^{k(p-2)_-}\int_{2R>|y-z|>a_k}\frac{dy
}
{|y-z|^{d+\alpha-1}} dy\\
& \le c_2(R) b^{-k}2^{3k(p-2)_-+1}\le c_2 2^{-k((p-2)_{-}+2)},
\end{align*}
  where in the first inequality we used the mean value theorem  for the difference  inside the integral in the numerator, in the second last inequality, the fact $a_k\ge 2^{-k}$ and  in the last inequality, the fact $b\ge 2^{4(p-2)_- +3}$.\bk

For $II$, we have
$$
|II| \le c_3 \int_{D(R+b^{-k}, R+b^{-k}), y_d<5^{-k}}\frac{y_d^p+z_d^p+b^{-kp}}{|y-z|^{d+\alpha}}\sB(y, z)dy.
$$
Similarly as in the proof of Lemma \ref{l:Lvk}  we first estimate 
\begin{align*}
&\int_{D(R+b^{-k}, R+b^{-k}), y_d<5^{-k}}\frac{y_d^p+b^{-kp}}{|y-z|^{d+\alpha}}\sB(y, z)dy\\
&\le c_4\int_{\R^{d-1}}\int^{5^{-k}}_0\frac{y_d^{p+\beta_1/2}+10^{-kp}y_d^{\beta_1/2}}{|y-z|^{d+\alpha+\frac12\beta_1}}dy_dd{\widetilde y}\\
&\le c_5 5^{-k(p+1+\frac12\beta_1)}\int^\infty_0\frac{t^{d-2}}{(t^{2}+cz_d^2)^{(d+\alpha+\frac12\beta_1)/2}}dt \\
&= c_6 5^{-k(p+1+\frac12\beta_1)}z_d^{-1-\alpha-\frac12\beta_1}\int^\infty_0\frac{s^{d-2}}{(s^{2}+1)^{(d+\alpha+\frac12\beta_1)/2}}ds\\
& \le c_7 (2/5)^{k(p+1+\frac12\beta_1)}z_d^{\beta_1}.
\end{align*}
For the remaining part, we use a similar argument:
\begin{align*}
&\int_{D(R+b^{-k}, R+b^{-k}), y_d<5^{-k}}\frac{z_d^p}{|y-z|^{d+\alpha}}\sB(y, z)dy\\
&\le c_8 z_d^p\int_{\R^{d-1}}\int^{5^{-k}}_0\frac{y_d^{\beta_1/2}}{|y-z|^{d+\alpha+\frac12\beta_1}}dy_dd\widetilde{y}\\
&\le c_9 z_d^p 5^{-k(1+\frac12\beta_1)}\int^\infty_0\frac{t^{d-2}}{(t^2+cz^2_d)^{(d+\alpha+\frac12\beta_1)/2}}dt\\
&=c_{10} z_d^p 5^{-k(1+\frac12\beta_1)}z_d^{-1-\alpha-\frac12\beta_1}\int^\infty_0\frac{s^{d-2}}{(s^2+1)^{(d+\alpha+\frac12\beta_1)/2}}ds \\
&\le c_{11} (4/5)^{k(1+\frac12\beta_1)}2^{-k(2+\beta_1)}z_d^{p-1-\alpha-\frac12\beta_1}\le c_{12} (4/5)^{k(1+\frac12\beta_1)}z_d^{\beta_1}.
\end{align*}
Thus
$$
|II|\le c_{13} \big((2/5)^{k(p+1+\frac12 \beta_1)}+(4/5)^{k(1+\frac12\beta_1)}\big)z_d^{\beta_1}.
$$

Let $B_k:=U(R)\cap \{y_d>5^{-k}\}\cap \{y:|y-z|>2^{-1}z_d\}$. Then, we have
\begin{align*}
|III|&=\left|\int_{B_k}\frac{(y_d-u_d)^p-y_d^p-((z_d-u_d)^p-z_d^p)}{|y-z|^{d+\alpha}}
\sB(y, z)dy\right|\\
&=\left|\int_{B_k}\frac{
 F_1(y_d, z_d, u_d) \bk}{|y-z|^{d+\alpha}}
\sB(y, z)dy\right|\\
&\le c_{14}b^{-k}2^{k3(p-2)_-}\int_{U(R), |y-z|>2^{-1}z_d}\frac{1}{|y-z|^{d+\alpha-1}}dy\\
&\le c_{15}(4/5)^{k}2^{-3k}(z_d^{1-\alpha}\vee \log \frac1{z_d})
\le c_{16} (4/5)^{k}z_d^3(z_d^{1-\alpha}\vee \log \frac1{z_d}),
\end{align*}
where in the first inequality we use the mean value theorem inside the integral in the numerator
and the fact the derivative of the integrand is bounded above by $c(5^{-k}-b^{-k})^{-(p-2)_-}\le c2^{3k(p-2)_-}$.

Let $F_2(y_d, z_d, u_d):=p(p-1)\int^1_0\left((z_d-u_d+t(y_d-z_d))^{p-2}-(z_d+t(y_d-z_d))^{p-2}\right)(1-t)dt$. 
For $t\in [0, 1]$, $u\in B(0, b^{-k})$
and $y\in B(z, \frac12 z_d)$, $z_d-u_d+t(y_d-z_d)$ and $z_d+t(y_d-z_d)$ are both comparable with $z_d$. Thus,
for $IV$, we have that,  for large $k$,

\begin{align*}
&|IV| \le \left|
 \mathrm{p.v.} 
\int_{U(R)\cap B(z, 2^{-1}z_d)}\frac{(y_d-u_d)^p-y_d^p-((z_d-u_d)^p-z_d^p)}{|y-z|^{d+\alpha}}
(\sB(y, z)-1)dy\right|\\
&+\left|\int_{U(R)\cap B(z, 2^{-1}z_d)}
\left[(y_d-u_d)^p-y_d^p-((z_d-u_d)^p-z_d^p) 
\right.
 \right.\\
&\qquad\qquad\qquad \left.\left.
-p((z_d-u_d)^{p-1}-z_d^{p-1})
{\bf 1}_{B(z, a_k)} (y) (y_d-z_d) \right]{|y-z|^{-d-\alpha}}
dy\right|\\
&\le 
\int_{U(R)\cap B(z, 2^{-1}z_d)}
\frac{\left|F_1(y_d, z_d, u_d)\right|}{z_d^{ \theta} |y-z|^{d+\alpha- \theta }}
+
\int_{U(R)\cap (B(z, 2^{-1}z_d) \setminus B(z, a_k) )}
\frac{ \left|F_1(y_d, z_d, u_d)\right|}{ |y-z|^{d+\alpha}}
dy\\
& \quad +\int_{B(z, a_k)}\frac{\left| F_2(y_d, z_d, u_d) \right|}{|y-z|^{d+\alpha-2}}
dy\\
&\le c_{17} \left(
\int_{B(z, 2^{-1}z_d)}
\frac{ |u_d|z_d^{p-2}
dy}
{z_d^\theta |y-z|^{d+\alpha-1-\theta}}
+ 
\int_{B(z, 2^{-1}z_d) \setminus B(z, a_k) }
\frac{
 |u_d|z_d^{p-2}dy
}{ |y-z|^{d+\alpha-1}}
+\int_{B(z, a_k)}\frac{ |u_d|z_d^{p-3}dy}{|y-z|^{d+\alpha-2}}
 \right)\\
&\le c_{ 18 }b^{-k} \left ( z_d^{p-1-\alpha} +z_d^{p-2 } 
a_k^{1-\alpha} +
z_d^{p-1-\alpha}  \right)   
\le  c_{19} (2/b)^k  z_d^{p-1-\alpha }
\le c_{ 20 } (4/5)^k z_d^{p+1-\alpha },
\end{align*}
where  in the second inequality we used 
\textbf{(A2)}  and   in the third inequality 
the mean value theorem
 for the difference inside the integral in the definitions of $F_1$ and $F_2$.

 Combining the estimates for $I$, $II$, $III$ and $IV$, we arrive at the desired assertion. \qed

\begin{lemma}\label{l:estimate-of-L-hat-B}
\begin{itemize}
\item [(a)]
  There exists $C_1>0$ such that for every $R \ge 1$ and $z\in U(R)$, 
$$
0\ge L_{\alpha}^\sB  h_{\alpha-1, R}(z) \ge -C_1 z_d^{\beta_1}(|\log z_d|^{\beta_3}\vee1)  \int_{|y|\ge 
R \normal
}|y|^{-\beta_1-d-1}\big(1+{\bf 1}_{|y|\ge1}(\log|y|)^{\beta_3}\big)\, dy\, .
$$
\item [(b)] Let $\alpha-1 <p <\alpha+\beta_1$. There exist $r_0\in (0,1/2]$ and  $C_2>0$ and $C_3>0$ such that for every  $z\in D(\frac12, r_0)$,
$$
C_{2} z_d^{p-\alpha}\le L_{\alpha}^\sB  h_p(z) \le C_{3} z_d^{p-\alpha}.
$$
\end{itemize}
\end{lemma}

\pf (a) Let $z\in U (R)$. 
Then by \eqref{e:Lal},
we see that
\begin{align*}
L_{\alpha}^\sB  h_{\alpha-1, R}(z)=-\int_{D(R,R)^c\cap \R^d_+}\frac{y_d^{ \alpha -1} }{|y-z|^{d+\alpha}}\sB (z,y)\, dy\, ,
\end{align*}
which is negative. Further, if $y\in D(R,R)^c$ and $z\in U(R)=D(R/2,R/2)$, then 
$|y|\ge   R>|z|$, $|y-z|\ge z_d$ and $|y-z| \asymp |y|$. 
 Thus it follows from 
\cite[Lemma 5.2]{KSV} that
\begin{align*}
&|L_{\alpha}^\sB  h_{\alpha-1, R}(z)|\le 
\int_{D(R,R)^c\cap \R^d_+}
\frac{|y|^{\alpha-1}}{|y-z|^{d+\alpha}}\sB (z,y)\, dy \\
&\le c_1\int_{y \in \R^d_+, |y|\ge  R, |y-z|\ge z_d}
|y|^{-d-1}\sB (z,y)\, dy \\
&\le c_2 z_d^{\beta_1}(|\log z_d|^{\beta_3}\vee1)  
\int_{|y|\ge 
 R}
|y|^{-\beta_1-d-1}\big(1+{\bf 1}_{|y|\ge1}(\log|y|)^{\beta_3}\big)\, dy
\, .
\end{align*}

\noindent (b) Let $z\in U$. Then by using Lemma \ref{l:LB-on-g}, 
\begin{eqnarray*}
L_{\alpha}^\sB  h_p(z)
=C(\alpha, p, \sB )z_d^{p-\alpha}-\int_{D(1,1)^c\cap \R^d_+}\frac{y_d^p}{|y-z|^{d+\alpha}}\sB (z, y)\, dy\, .
\end{eqnarray*}
Since the second term is non-negative, by removing it we obtain the upper bound.
In the same way as before (this uses $p<\alpha+\beta_1$),
$$
\left| \int_{D(1,1)^c\cap \R^d_+}\frac{y_d^p}{|y-z|^{d+\alpha}}\sB (z, y)\, dy \right| \le c_2 z_d^{\beta_1}|\log z_d|^{\beta_3}\, .
$$
Thus, for any $z\in U$,
$$
L_{\alpha}^\sB h_p(z) 
\ge  C(\alpha, p, \sB )z_d^{p-\alpha}-c_2 z_d^{\beta_1}|\log z_d|^{\beta_3}\, .
$$
Since $p-\alpha<\beta_1$ and $C(\alpha, p, \sB )>0$, we can find $r_0\in (0,1/2]$ such that the function $t\mapsto C(\alpha, p, \sB ))
t^{p-\alpha}-c_2 t^{\beta_1}|\log t|^{\beta_3} $ 
 is  
 bounded from below by $c_3 t^{p-\alpha}$ with a positive constant $c_3>0$ for all $t\in (0,r_0)$. 
This concludes the proof of the lower bound. 
\qed

In the remainder of this paper, $r_0$ always stands for the constant in the lemma above.

\begin{lemma}\label{l:exit-time estimate-U}
There exists a constant $C>0$  such that
\begin{equation}\label{e:exit-time-estimate-U}
\E_x \tau_U \le C
x_d^{\alpha-1}\,, \quad x\in U.
\end{equation}
\end{lemma}

\pf
Choose $q\in (\alpha-1, \alpha)$ and let $\eta(x):=h_{\alpha-1}(x)-h_q(x)$, $x\in \R^d_+$.
For $x\notin D(1,1)$, $\eta(x)=0$, while if $x\in D(1,1)$ we have $\eta(x)=x_d^{\alpha-1}-x_d^q>0$. By Lemma \ref{l:estimate-of-L-hat-B}, for all $x\in U(r_0)$ we have that 
$L^\sB h_{\alpha-1}(x)\le 0$ and $L^\sB h_q(x)\ge c_1 x^{q-\alpha}$. Thus we can find $r_1\in (0, r_0]$ such that
\begin{equation}\label{e:exit-time-estimate-U-3}
L_{\alpha}^\sB \eta(x)=L_{\alpha}^\sB h_{\alpha-1}(x)-L_{\alpha}^\sB h_q(x)\le -c_1 x_d^{q-\alpha}\le -1, \qquad x\in U(r_1).
\end{equation}
Let $g_k=v_k\ast (\eta {\bf 1}_{D(5^{-k}})$. It follows from Lemmas \ref{l:Lvk} and \ref{l:hq} applied to $h_q$ and $h_{\alpha-1}$ that  $L_{\alpha}^\sB g_k\to L_{\alpha}^\sB \eta$ on $U$ and the sequence of functions
$|L_{\alpha}^\sB g_k-L_{\alpha}^\sB \eta|$  is bounded by some constant $c_2>0$. In particular, 
\begin{equation}\label{e:exit-time-estimate-U-2}
-L_{\alpha}^\sB g_k(z) \ge -L_{\alpha}^\sB \eta(z)-c_2 \ge 1-c_2\, , \qquad z\in U(r), r\le r_1\, .
\end{equation}
It follows from \cite[Lemma 3.6]{KSV} that for all $t\ge 0$,
$$
\E_x [g_k(Y_{t\wedge \tau_{U_{a_k}(r)}})]-g_k(x)=\E_x\int_0^t {\bf 1}_{s<\tau_{U_{a_k}(r)}} L_{\alpha}^\sB g_k (Y_s)\, ds \, , \qquad x\in U(r), r\le r_1.
$$
As $k\to \infty$, the left-hand side converges to $\E_x [\eta(Y_{t\wedge \tau_{U(r)}})]-\eta(x)$. For the right-hand side we can use Fatou's lemma (justified because of \eqref{e:exit-time-estimate-U-2}) to conclude that for $x\in U(r)$ with $ r\le r_1$,
\begin{align*}
\limsup_{k\to \infty}\E_x\int_0^t {\bf 1}_{s<\tau_{U_{a_k}(r)}} L_{\alpha}^\sB g_k(Y_s)\, ds \le \E_x\int_0^t {\bf 1}_{s<\tau_{U(r)}} L_{\alpha}^\sB \eta(Y_s)\, ds \le -\E_x (t\wedge \tau_{U(r)})\, .
\end{align*}
Thus we get that $\E_x [\eta(Y_{t\wedge \tau_{U(r)}})]-\eta(x) \le -\E_x (t\wedge \tau_{U(r)})$, and by letting $t\to \infty$,
$$
-\eta(x)\le \E_x [\eta(Y_{\tau_{U(r)}})]-\eta(x) \le -\E_x \tau_{U(r)}\, \, , \qquad x\in U(r), r\le r_1.
$$
Thus we get $\E_x \tau_{U(r)}\le \eta(x)\le x_d^{\alpha-1}$. By using that $U(r_1)=r_1 U$ and \eqref{e:exit-time-scaling},
 for any $x\in U$,
$$
\E_x \tau_U=r_1^{-\alpha} \E_{r_1x}\tau_{r_1U}\le r_1^{-\alpha}(r_1x_d)^{\alpha-1}= r_1^{-1} x_d^{\alpha-1}\, .
$$
We have proved the claim of the lemma with $C=r_1^{-1} $. \qed
\normal

\begin{prop}\label{l:dynkin-hp}
Let $p\in [\alpha-1, \alpha+\beta_1)$, $R \ge 1$ and $r\le R$. For every $x\in U(r)$ it holds that
\begin{equation}\label{e:dynkin-hp}
\E_x [h_{p, R}(Y_{ \tau_{U(r)}})]=h_{p, R}(x)+\E_x \int_0^{\tau_{U(r)}}  L_{\alpha}^\sB h_{p, R}(Y_s)\, ds \, .
\end{equation}
\end{prop}
\pf 
Set $g_k:=v_k\ast (h_{p, R} {\bf 1}_{D(5^{-k})})$. Let $x\in U(r)$, 
$r\le  R$. 
There is $k_0\ge 1$ such that $x\in U_{a_k}(r)$ for all $k\ge k_0$. Note that since $g_k\in C_c^{\infty}(\R^d_+)$, it follows from \cite[Lemma 3.6]{KSV} that for all $t\ge 0$,
$$
\E_x [g_k(Y_{t\wedge \tau_{U_{a_k}(r)}})]=g_k(x)+\E_x\int_0^t {\bf 1}_{s<\tau_{U_{a_k}(r)}} L_{\alpha}^\sB g_k (Y_s)\, ds \, .
$$
Clearly, $\lim_{k\to \infty}\tau_{U_{a_k}(r)}=\tau_{U(r)}$. Since $g_k\to h_{p, R}$ as $k\to \infty$, we get that the left-hand side above converges to $\E_x [h_{p, R}(Y_{t\wedge \tau_{U(r)}})]$.

On the other hand, 
by combining Lemmas \ref{l:Lvk} and \ref{l:hq}, we see that for every $z\in U(r)$ with $r\le R$,
$$
\lim_{k\to \infty}|L_{\alpha}^\sB g_k (z)- L_{\alpha}^\sB h_{p, R}(z)|=0
$$
and $|L_{\alpha}^\sB g_k (z)- L_{\alpha}^\sB h_{p, R}(z)|$
is bounded. 
Thus, we can use the bounded convergence theorem and get  
$$
\lim_{k\to \infty} \E_x\int_0^t {\bf 1}_{s<\tau_{U_{a_k}(r)}} L_{\alpha}^\sB g_k (Y_s)\, ds = \E_x \int_0^t {\bf 1}_{s<\tau_{U(r)}} L_{\alpha}^\sB h_{p, R}(Y_s)\, ds \, .
$$
Therefore, 
$$
\E_x [h_{p, R}(Y_{t\wedge \tau_{U_{a_k}(r)}})]=h_{p, R}(x)+\E_x\int_0^{t \wedge \tau_{U(r)}} L_{\alpha}^\sB h_{p, R} (Y_s)\, ds \, .
$$
By letting $t\to \infty$ we get that the left-hand side above converges to $\E_x [h_{p, R}(Y_{\tau_{U(r)}})]$.

When $p=\alpha-1$, by Lemma \ref{l:estimate-of-L-hat-B} (a), $ L_{\alpha}^\sB h_{p, R}(z)\le 0$. Thus
we can use use the monotone convergence theorem and obtain \eqref{e:dynkin-hp}.

When 
$p\in (\alpha-1, \alpha+\beta_1)$, by Lemma \ref{l:estimate-of-L-hat-B} (b) and 
scaling, we have 
$ L_{\alpha}^\sB h_{p, R}(z)> 0$
on $D(R/2, r_0R)\supset D(r/2, rr_0/2)$.
Thus, we can use the monotone convergence theorem and get
$$
\lim_{t\to \infty}\E_x\int_0^{t \wedge \tau_{U(r)}} {\bf 1}_{Y_s \in D(r/2, rr_0/2)} 
L_{\alpha}^\sB h_{p, R} (Y_s) \, ds 
=\E_x \int_0^{\tau_{U(r)}}  {\bf 1}_{Y_s \in D(r/2, rr_0/2)} L_{\alpha}^\sB h_{p, R}(Y_s)\, ds.
$$
On the other hand,
since
$$
L_{\alpha}^\sB h_{p, R}(z)=C(\alpha, p, \sB)z_d^{p-\alpha}-\int_{D(R, R)}\frac{y^p_d}{|y-z|^{d+\alpha}}\sB(y, z)dy, \quad z\in D(R, R),
$$
we know that
 $L_{\alpha}^\sB h_{p, R}(z)$ is bounded on $U(r) \setminus D(r/2, rr_0/2)$.
Thus, using Lemma \ref{l:exit-time estimate-U} and the bounded convergence theorem, we get
$$
\lim_{t\to \infty}\E_x\int_0^{t \wedge \tau_{U(r)}} {\bf 1}_{Y_s \in U(r) \setminus D(r/2, rr_0/2)}
L_{\alpha}^\sB h_{p, R} (Y_s)\, ds 
=\E_x \int_0^{\tau_{U(r)}}  {\bf 1}_{Y_s \in U(r) \setminus D(r/2, rr_0/2)} L_{\alpha}^\sB h_{p, R}(Y_s)\, ds.
$$
Combining these, we obtain \eqref{e:dynkin-hp} for $p\in (\alpha-1, \alpha+\beta_1)$ too.
\qed

\section{Boundary Harnack principle} \label{s:BHP}

In this section we always assume that $\alpha>1$.
The next two results are applications of Proposition  \ref{l:dynkin-hp}.

\begin{lemma}\label{l:gharmonic}
 For all $r>0$ it holds that 
\begin{equation}\label{e:gharmonic}
\E_x [h_{\alpha-1, \infty}(Y_{ \tau_{U(r)}})]=h_{\alpha-1, \infty}(x), \quad \text{for all } x \in U(r).
\end{equation}
 In particular, the function $h_{\alpha-1, \infty}(x)=g_{\alpha-1}(x)=x_d^{\alpha-1}$ is harmonic  in $\R^d_+$ with respect to $Y$.  
\end{lemma}
\pf
Fix  $r>0$ and $x \in U(r)$. Let  $R \ge 1 \vee r$.
By Proposition  \ref{l:dynkin-hp},
\begin{equation}\label{e:ha}
\E_x [h_{\alpha-1, R}(Y_{ \tau_{U(r)}})]=x_d^{\alpha-1}+\E_x \int_0^{\tau_{U(r)}}  L_{\alpha}^\sB h_{\alpha-1, R}(Y_s)\, ds \, .
\end{equation}
By the monotone convergence theorem, 
$$\lim_{R \to \infty}
\E_x [h_{\alpha-1, R}(Y_{ \tau_{U(r)}})]=\E_x [h_{\alpha-1, \infty}(Y_{ \tau_{U(r)}})].$$
Using 
Lemmas \ref{l:estimate-of-L-hat-B} (a) and \ref{l:exit-time estimate-U}, we see that 
\begin{align*}
0 &\ge \E_x \int_0^{\tau_{U(r)}}  L_{\alpha}^\sB h_{\alpha-1, R}(Y_s)\, ds\\
 &\ge -c_1(r)\E_x {\tau_{U(r)}} \int_{|y|\ge 
 R}
 |y|^{-\beta_1-d-1}\big(1+{\bf 1}_{|y|\ge1}(\log|y|)^{\beta_3}\big)\, dy\\
 &\ge -c_2(r) \int_{|y|\ge 
  R}
 |y|^{-\beta_1-d-1}\big(1+{\bf 1}_{|y|\ge1}(\log|y|)^{\beta_3}\big)\, dy.
\end{align*}
Thus 
$$
\lim_{R \to \infty}\E_x \int_0^{\tau_{U(r)}}  L_{\alpha}^\sB h_{\alpha-1, R}(Y_s)\, ds=0.
$$
The proof is now complete.
\qed

Recall that  $r_0$ in the constant in Lemma \ref{l:estimate-of-L-hat-B}(b). 
\begin{lemma}\label{l:exltlow}
There exists $C>0$ such that 
$$
\E_x\tau_{U(r_0)} \le C
\P_x\left(Y_{\tau_{U(r_0)}}\in  D(1, 1)  \right)  \quad \text{for all } x \in U(r_0).
$$
\end{lemma}
\pf
Choose a  $p\in (\alpha-1, \alpha)$.
By Proposition \ref{l:dynkin-hp}, for every $x\in U(r_0)$,
\begin{align*}
\E_x [h_{p}(Y_{ \tau_{U(r_0)}})]=h_{p}(x)+\E_x \int_0^{\tau_{U(r_0)}}  L_{\alpha}^\sB h_{p}(Y_s)\, ds \, .
\end{align*}
Thus, using Lemma \ref{l:estimate-of-L-hat-B} (b) and that $h_{p}$ is bounded by 1 and supported 
on $D(1, 1)$, we get
$$
\P_x\left(Y_{\tau_{U(r_0)}}\in D(1, 1)\right)   \ge \E_x [h_{p}(Y_{ \tau_{U(r_0)}})]  \ge \E_x \int_0^{\tau_{U(r_0)}}  L_{\alpha}^\sB h_{p}(Y_s)\, ds 
\ge c \E_x\tau_{U(r_0)}.
$$
\qed

Since $\int_{D(1,1)}y_d^{\alpha-1} |y|^{-d-\alpha-\beta_1}\, dy =+\infty$, by the same argument as that leading to  \cite[(5.10)]{KSV}, we also have that 
there exists an 
$r_1 \in (0, r_0)$ small enough so that 
for all $r\in (0, r_1]$  and $ R \ge 1$, 
\begin{equation}\label{e:integral-upper-g} 
\E_x \int_0^{\tau_{U(r)}} (Y_t^d)^{\beta_1} |\log Y_t^d|^{\beta_3} \, dt \le \E_x[h_{\alpha-1}\normal(Y_{\tau_{U(r)}})] 
\le \E_x[h_{ \alpha-1, R}(Y_{\tau_{U(r)}})]\, ,\normal  \quad x\in U(r) \, .
\end{equation}

Using \eqref{e:exit-time-scaling},  \eqref{e:integral-upper-g}, Proposition \ref{l:dynkin-hp} and Lemma \ref{l:estimate-of-L-hat-B} (a), 
repeating the proof of \cite[Lemma 5.7]{KSV}, we get
the following upper bound.

\begin{lemma}\label{l:upper-bound-for-integral}
There exists a constant $C>0$ such that for
all $x\in U$,
\begin{equation}\label{e:upper-bound-for-integral}
\E_x \int_0^{\tau_{U}} (Y_t^d)^{\beta_1}  |\log Y_t^d|^{\beta_3}\, dt \le C x_d^{\alpha-1}\, .
\end{equation}
\end{lemma}

Using the above lemma, \cite[Lemma 5.2(b)]{KSV} and the scaling 
property of $Y$, the proof of the next lemma is the same as that of \cite[Lemma 6.1)]{KSV}, we omit the proof.
\begin{lemma}\label{e:POTAe7.14}
There exists $C>0$ such that for all 
$0<4r\le R\le 1$ and $w\in D(r,r)$, 
$$
\P_w\Big(Y_{\tau_{B(w, r)\cap \R^d_+}}\in A(w, R, 4)\cap \R^d_+\Big)\le C
\frac{r^{\alpha+\beta_1}}{R^{\alpha+\beta_1}}\frac{w_d^{\alpha-1}}{r^{\alpha-1}}.
$$
\end{lemma}

Recall that  $r_0$ in the constant in Lemma \ref{l:estimate-of-L-hat-B}(b). 
\begin{lemma}\label{l:POTAl7.4}
There exists $C>0$  such that 
for any $x \in U(2^{-4}r_0)$,
$$
\P_x\left(Y_{\tau_{U(r_0)}}\in D(1, 1)\right)\le C
\P_x\left(Y_{\tau_{U(r_0)}}\in D(1/2, 1)\setminus D(1/2, 3/4)\right).
$$
\end{lemma}

\pf 
Let $V=U(r_0)$ and 
$$
H_2:=\{Y_{\tau_{U}}\in D(1, 1)\}, \quad H_1:=\{Y_{\tau_{U}}\in D(1/2, 1)\setminus D(1/2, 3/4)\}.
$$

Choose a $p \in (\alpha-1, \alpha)$ and let
$
\kappa(x)
=C(\alpha, p,  \sB)x_d^{-\alpha}.
$
Let $Y^{\kappa}$ be the subprocess of $Y$ with killing potential $\kappa$ so that the corresponding Dirichlet form  is
$
\EE(u,v)+\int_{\R^d_+} u(x)v(x)\kappa(x) dx .
$
Either by repeating the proof of \cite[(5.17)]{KSV} or using \cite[Theorem 1.3]{KSV}, we get that
$$
\P_w\left(Y^{\kappa}_{\tau_{V}}\in D(r_0/4, 1)\setminus D(r_0/4, 3/4)\right)
 \ge c_1 w^p_d, \quad w\in U(r_0/2).
$$
Thus,
\begin{equation}\label{e:POTAe7.15}
\P_w(H_1)\ge 
\P_w\left(Y^{\kappa}_{\tau_{V}}\in D(1/2, 1)\setminus D(1/2, 3/4)\right)
 \ge c_1 w^p_d, \quad 
w\in U(r_0/2).
\end{equation}
For $i\ge 1$, set
$$
s_0=s_1, \quad s_i=\frac{r_0}8\Big(\frac12-\frac1{50}\sum^i_{j=1}\frac1{j^2}\Big)\quad \text{ and }  \quad J_i=D(s_i, 2^{-i-3}r_0)\setminus D(s_i, 2^{-i-4}r_0) .
$$
Note that $r_0/(20)<s_i<r_0/(16)$. Define for $i\ge 1$,
\begin{equation}\label{e:POTAe7.16}
d_i=\sup_{z\in J_i}\frac{\P_z(H_2)}{\P_z(H_1)}, 
\quad \widetilde{J}_i=
D(s_{i-1}, 2^{-i-3}r_0),
\quad \tau_i=\tau_{\widetilde{J}_i}.
\end{equation}
Repeating the argument leading to \cite[(6.29)]{KSV19}, we get that for 
$z\in J_i$ and $i\ge 2$,
\begin{equation}\label{e:POTAe7.17}
\P_z(H_2)\le \Big(\sup_{1\le k\le i-1}d_k\Big)\P_z(H_1)
+\P_z\left(Y_{\tau_i}\in 
D(1, 1)
\setminus \cup^{i-1}_{k=1}J_k\right).
\end{equation}

Recall that  $n_0$ in the constant in Corollary \ref{c:lemma4-1} (b). 
For  $i\ge 2$, define $\sigma_{i,0}=0, \sigma_{i,1}=\inf\{t>0: |Y_t-Y_0|\geq n_0 2^{-i-2}r_0\} $ and
$\sigma_{i,m+1}=\sigma_{i,m}+\sigma_{i,1}\circ\theta_{\sigma_{i,m}}$
for $m\geq 1$. By Corollary \ref{c:lemma4-1} (b),  we have that 
\begin{align}\label{e:POTAe7.18}
\P_{w}(Y_{\sigma_{i,1}}\in \widetilde{J}_i) \le 1-
\P_{w}( \sigma_{i,1}=\zeta ) \le 1-\P_{w}(
\tau_{B(w,n_0 w_d)}
=\zeta )  <2^{-1},\ \ \ 
w\in  \widetilde{J}_i.
\end{align}
For the purpose of further estimates, we now choose a
positive integer $l$ 
such that 
$l \ge \alpha+\beta_1$.
Next we choose 
$i_0 \ge 2$ large enough so that
$n_02^{-i+1}<1/(200 l i^3)$ for all $i\ge i_0$.
Now we assume $i\ge i_0$.
Using \eqref{e:POTAe7.18} and the strong Markov property we have
that  for   $z\in J_i$,
\begin{align}\label{e:POTAe7.19}    
&\P_z( \tau_{i}>\sigma_{i,li})\leq \P_z(Y_{\sigma_{i,k}}\in \widetilde{J}_i, 1\leq k\leq
li
 )\nonumber\\
 &=
 \E_z \left[ \P_{Y_{\sigma_{i,li-1}}} (Y_{\sigma_{i,1}}\in  \widetilde{J}_i) : Y_{\sigma_{i,li-1}} \in  \widetilde{J}_i,  Y_{\sigma_{i,k}}\in \widetilde{J}_i, 1\leq k\leq
li-2
  \right]
 \nonumber\\      
 &\leq \P_z\left(Y_{\sigma_{i,k}}\in \widetilde{J}_i, 1\leq k\leq
li-1
 \right)2^{-1}\leq 2^{-li}.
\end{align}
Note that if $z\in J_i$ and $y\in 
D(1, 1) 
\setminus[ \widetilde{J}_i \cup(\cup_{k=1}^{i-1}J_k)]$, 
then $|y-z|\ge (s_{i-1}-s_i) \wedge (2^{-4}r_0-2^{-i-3}r_0)  = r_0/(400 i^2)$.
Furthermore, 
since  $2^{-i-2}n_0 < 1/(400 i^2)$ (recall that $i\ge  i_0$), 
if $Y_{\tau_i}(\omega)\in 
D(1, 1)
\setminus \cup_{k=1}^{i-1}J_k$ and $\tau_i(\omega)\le \sigma_{i,li}(\omega)$, then $\tau_i(\omega)=\sigma_{i,k}(\omega)$ for some $k=k(\omega)\le li$.  
Dependence of $k$ on $\omega$ will be omitted  in the next few lines. 
Hence on
$\{Y_{\tau_{i}}\in 
D (1 , 1) 
\setminus \cup_{k=1}^{i-1}J_k,\ \ \tau_{i}\leq
\sigma_{i,li}\}$ with $Y_0=z\in J_i$,  we have
$|Y_{\sigma_{i,k}}-Y_{\sigma_{i,0}}|=|Y_{\tau_i}-Y_0|>
\frac{r_0}{400i^2}$ for some $1\leq k\leq li$.
Thus  for some $1\leq k\leq li$,
$
\sum_{j=0}^k|Y_{\sigma_{i,j}}-Y_{\sigma_{i,j-1}}|>
r_0(400i^2)^{-1}
$
which implies for some $1\leq j\leq k\le  li$,
$ |Y_{\sigma_{i,j}}-Y_{\sigma_{i,j-1}}|\geq
r_0({k400i^2})^{-1}\ge r_0(li)^{-1} (400i^2)^{-1}\, .
$
Thus, we have 
\begin{align*}
& \{Y_{\tau_{i}}\in 
D(1, 1) 
\setminus
\cup_{k=1}^{i-1}J_k,\ \ \tau_{i}\leq \sigma_{i,li}\}\\
  \subset  &  \cup_{j=1}^{li}\{|Y_{\sigma_{i,j}}- Y_{\sigma_{i,j-1}}|\geq r_0/(800li^3),
  Y_{\sigma_{i,j}}\in 
 D(1, 1), 
  Y_{\sigma_{i,j-1}}\in \widetilde{J}_{i}
\}.
\end{align*}
Now, using  Lemma \ref{e:POTAe7.14} (with $r=2^{-i-2}n_0r_0$ and $R=r_0/(800 l i^3)$) 
 (noting that  $4 \cdot 2^{-i-1}n_0<1/(400 l i^3)$ for all $i\ge i_0$), 
 and repeating the argument leading to \cite[(6.34)]{KSV19}, 
 we get that for $z\in J_i$,
\begin{align*}
&\P_z \left(Y_{\tau_{i}}\in D(1, 1) \setminus
\cup_{k=1}^{i-1}J_k,\ \ \tau_{i}\leq \sigma_{i,li} \right) \leq  li \sup_{w\in \widetilde{J}_{i}} 
  \P_w\left(|Y_{\sigma_{i,1}}-w  |\geq r_0(800li^3)^{-1},  Y_{\sigma_{i,1}}\in 
D(1, 1)
\right)\\
&\le  li \sup_{w\in \widetilde{J}_{i}} 
  \P_w\left(4>|Y_{\sigma_{i,1}}-w  |\geq r_0(800li^3)^{-1}\right)
 \le c_{12}li  \left(
\frac{800li^3}{2^{i+3}}\right)^{\alpha+\beta_1}.
\end{align*}
By this and (\ref{e:POTAe7.19}),  we have for 
$z\in J_i$,  $i\ge i_0$,
\begin{align}\label{e:POTAe7.22}
&\P_z\left( Y_{\tau_{i}}\in 
D(1, 1) 
\setminus \cup_{k=1}^{i-1}J_k \right) \leq 2^{-li} 
+c_{2} li  \left(
\frac{800li^3 }{2^{i+3}}\right)^{\alpha+\beta_1}. 
\end{align} 
By our choice of $l$, we have 
\begin{align}\label{e:POTAe7.23}
&li    \left(
\frac{800li^3}{2^{i+3}}\right)^{\alpha+\beta_1} 
= 100^{\alpha+\beta_1} l^{1+\alpha+\beta_1} i^{1+3(\alpha+\beta_1)} \left(2^{-(\alpha+\beta_1)}\right)^i
\ge \left(2^{-(\alpha+\beta_1)}\right)^i
 \ge (2^{-l})^{i}. 
\end{align}
Thus combining \eqref{e:POTAe7.23} with \eqref{e:POTAe7.22}, and then 
 using \eqref{e:POTAe7.15},  we get that  for
$z\in J_i$, $i\ge i_0$, 
\begin{align}\label{e:POTAe7.24}
&\frac{\P_z( Y_{\tau_i}\in
D(1, 1) 
\setminus \cup_{k=1}^{i-1}J_k)}{\P_z(H_1)} \le 
c_{3} li  
2^{ip}
\left(
\frac{800li^3 }{2^{i+3}}\right)^{\alpha+\beta_1}
\le  
c_{4} i ^{1+3(\alpha+\beta_1)}2^{(p-\alpha-\beta_1)i}.
\end{align} 
By this and (\ref{e:POTAe7.17}),  for 
$z\in J_i$,  $i\ge i_0$,for all $i\ge i_0$
\begin{align*}  
&\frac{\P_z( H_2)}{\P_z(H_1)} \leq  \sup_{1\leq k\leq i-1}  d_k   +\frac{\P_z( Y_{\tau_i}\in
D(1 , 1) 
\setminus \cup_{k=1}^{i-1}J_k)}{\P_z(H_1)}\le \sup_{1\leq k\leq i-1}d_k+
c_{4} \frac{i ^{1+3(\alpha+\beta_1)}}{2^{(\alpha+\beta_1-p)i}}.
\end{align*}
This implies that for all $i\ge1$
\begin{align} 
d_i& \leq  \sup_{1\leq k\leq i_0-1} d_k
+c_{4}\sum_{k=1}^i\frac{i ^{1+3(\alpha+\beta_1)}}{2^{(\alpha+\beta_1-p)i}}
\leq
\sup_{1\leq k\leq i_0-1} d_k
+c_{4}\sum_{k=1}^\infty\frac{i ^{1+3(\alpha+\beta_1)}}{2^{(\alpha+\beta_1-p)i}}
=:c_{5} <\infty.\nonumber
\end{align}
Since $U(2^{-4}r_0) 
\subset \cup_{k=1}^\infty J_k$, the proof is now complete.
\qed

Using Corollary \ref{c:lemma4-1} (b), we can prove the following Carleson estimate.
\begin{thm}[Carleson estimate]\label{t:carleson}
There exists a constant $C>0$  
such that for any $w \in\partial \R^d_+$, $r>0$, 
and any non-negative  function $f$ in $\R^d_+$ that is harmonic in $\R^d_+ \cap B(w, r)$ with respect to $Y$ and vanishes continuously on $ \partial \R^d_+ \cap B(w, r)$, we have
\begin{equation}\label{e:carleson}
f(x)\le C f(\wh{x}) \qquad \hbox{for all }  x\in \R^d_+\cap B(w,r/2),
\end{equation}
where $\wh{x}\in \R^d_+\cap B(w,r)$ with $\wh{x}_d\ge r/4$.
\end{thm}
\pf
Recall that  $n_0$ in the constant in Corollary \ref{c:lemma4-1} (b).
By using $B_0(x)=B(x,n_0x_d)$ instead of $B_0(x)=B(x,x_d/2)$ in the proof of \cite[Theorem1.2]{KSV} and applying our 
 Corollary \ref{c:lemma4-1} (b), the proof of the theorem is almost identical to that of \cite[Theorem1.2]{KSV}. We omit the details. 
\qed
\normal

\noindent
{\bf Proof of Theorem \ref{t:BHP}.}
Recall that  $r_0$ in the constant in Lemma \ref{l:estimate-of-L-hat-B}(b).
 By scaling, it suffices to deal with the case $r=1$. Moreover, by
Theorem  \ref{t:uhp} (b), 
it suffices to prove \eqref{e:TAMSe1.8new}
for  $x, y\in D_{\wt w}(2^{-8}r_0, 2^{-8}r_0)$.
Since $f$ is harmonic in $D_{\wt w}(2, 2)$ and vanishes continuously on $B(\wt w, 2)\cap \partial \R^d_+$,
it is regular harmonic in 
$D_{\wt w}(7/4, 7/4)$ and vanishes continuously on 
$B(\wt w, 7/4)\cap \partial \R^d_+$ (see \cite[Lemma 5.1]{KSV19} and its proof).
Throughout the remainder of this proof, we assume that 
$x\in D_{\wt w}(2^{-8}r_0, 2^{-8}r_0)$. 
Without loss of generality we take $\widetilde{w}=0$.

Define $x_0=(\widetilde{x}, 1/(16))$ 
and $V=U(r_0)$. By the Harnack inequality 
and Lemma \ref{l:POTAl7.4}, we have
\begin{align}\label{e:TAMSe6.37}
f(x)&=\E_x[f(Y_{\tau_{V}})]\ge \E_x[f(Y_{\tau_{V}}); Y_{\tau_{V}}\in 
D(1/2, 1)\setminus D(1/2, 3/4)]
\nonumber\\
&\ge c_0f(x_0)\P_x(Y_{\tau_{V}}\in 
D(1/2, 1)\setminus D(1/2, 3/4))
\ge c_1f(x_0)\P_x(Y_{\tau_{V}}\in 
D(1, 1)).
\end{align}

Set $w_0=(\widetilde{0}, 2^{-7})$. 
Then, using \cite[Proposition 3.11 (a)]{KSV}, we also have  
\cite[(6.13)-(6.14)]{KSV}, that is, 
\begin{align}\label{e:POTAe7.27}
f(w_0)&\ge c_2 \int_{\R^d_+\setminus 
D(1, 1)}
 J \bk
({w_0}, y)f(y)dy,
\end{align}
and 
\begin{equation}\label{e:new-estimate-for-J}
 J \bk(z,y)\le c_{3}
 J \bk(w_0,y), \quad \text{ for any $z\in U$ and $y\in \R^d_+\setminus 
D(1, 1)$}.
\end{equation}

Combining \eqref{e:new-estimate-for-J} with  \eqref{e:POTAe7.27} we now have
\begin{align}\label{e:POTAe7.29}
&\E_x\left[f(Y_{\tau_{V}}); Y_{\tau_{V}}\notin 
D(1, 1)
\right]=\E_x\int^{\tau_{V}}_0\int_{\R^d_+\setminus 
D(1, 1)}
 J \bk(Y_t, y)f(y)dydt\nonumber\\
&\le c_{3} \E_x\tau_{V}\int_{\R^d_+\setminus 
D(1, 1)}
 J \bk(w_0, y)f(y)dy\le c_{4} f(w_0)\E_x\tau_{V}.
\end{align}

On the other hand, by the Harnack inequality (Theorem \ref{t:uhp}) and Carleson's estimate (Theorem \ref{t:carleson}), we have
\begin{align}\label{e:POTAe7.30}
\E_x\left[f(Y_{\tau_{V}}); Y_{\tau_{V}}\in 
D(1, 1)
\right]&\le c_{5}f(x_0)\P_x\left(Y_{\tau_{V}}\in 
D(1, 1)
\right).
\end{align}
Combining \eqref{e:POTAe7.29} and \eqref{e:POTAe7.30}, and using the Harnack inequality, we get
\begin{align*}
f(x)&=\E_x\left[f(Y_{\tau_{V}}); Y_{\tau_{V}}\in 
D(1, 1)
\right]+
\E_x\left[f(Y_{\tau_{V}}); Y_{\tau_{V}}\notin 
D(1, 1)
\right]\nonumber\\
&\le c_4f(w_0)\E_x\tau_{V}+c_{5}f(x_0)\P_x\left(Y_{\tau_{V}}\in 
D(1, 1)
\right) \le c_{6} f(x_0)\left(\E_x\tau_{V}+\P_x\left(Y_{\tau_{V}}\in 
D(1, 1)
\right) \right).
\end{align*}
 This with  
\eqref{e:POTAe7.27} and Lemma \ref{l:exltlow} implies that 
$
f(x) \asymp   f(x_0) \P_x\left(Y_{\tau_{V}}\in 
D(1, 1)
\right).
$
For any $y\in D(2^{-8}r_0, 2^{-8}r_0)$, 
we have the same estimate with $f(y_0)$ instead of $f(x_0)$, 
where $y_0=(\widetilde{y}, 1/(16))$. 
By the Harnack inequality, we have $f(x_0)\asymp f(y_0)$. 
Thus,
$$
\frac{f(x)}{f(y)}\asymp \frac{ \P_x\left(Y_{\tau_{V}}\in 
D(1, 1)
\right)} { \P_y\left(Y_{\tau_{V}}\in 
D(1, 1)
\right)}.
$$
We now apply this with $g_{\alpha-1}$ (which is harmonic by Lemma \ref{l:gharmonic}   ) to  conclude that 
$$
\frac{x_d^{\alpha-1}}{y_d^{\alpha-1}}\asymp \frac{ \P_x\left(Y_{\tau_{V}}\in 
D(1, 1)
\right)} { \P_y\left(Y_{\tau_{V}}\in 
D(1, 1)
\right)}\asymp  \frac{f(x)}{f(y)}.
$$
\qed

We now show that any non-negative function which is regular harmonic near a portion of boundary vanishes continuously on that portion of boundary,
cf.~\cite[Remark 6.2]{BBC} and \cite[Lemma 3.2]{CK02}. Thus, the above boundary Harnack principle also holds for regular harmonic functions.

\begin{lemma}\label{l:bh-decay}
There exists a constant $C>0$ such that for every bounded function  $f:\R^d_+\to [0,\infty)$  which is regular harmonic in $U=D(1/2,1/2)$, it holds that
$$
f(x)\le C \|f\|_{\infty} x_d^{\alpha-1}, \quad x\in D(2^{-5}, 2^{-5}).
$$
\end{lemma}
\pf Let  $f:\R^d_+\to [0,\infty)$ be a bounded function which is regular harmonic in $U=D(1/2,1/2)$. 
Then for every $x\in D(2^{-5}, 2^{-5})$, 
\begin{equation}\label{e:bh-decay-1}
f(x)=\E_x [f(Y_{\tau_U}), Y_{\tau_U}\in D(1,1)]+\E_x [f(Y_{\tau_U}), Y_{\tau_U}\notin D(1,1)].
\end{equation}
In the first term we use $f(Y_{\tau_U})\le \|f\|_{\infty}$ and then apply
Theorem \ref{t:BHP}
to get
\begin{equation}\label{e:bh.decay-2}
\E_x [f(Y_{\tau_U}), Y_{\tau_U}\in D(1,1)]\le
 \|f\|_{\infty} \bk
\P_x(Y_{\tau_U}\in D(1,1))\le c_1 \|f\|_{\infty}x_d^{\alpha-1}.
\end{equation}
Now we estimate the second term. For $z\in U$ and $w\in \R^d_+\setminus D(1, 1)$, we have $|w-z|\asymp |w|$. Thus, by using \cite[Lemma 5.2(a)]{KSV},
\begin{align*}
&\int_{\R^d_+\setminus D(1,1)}
 f(w) \bk
 \sB(z,w)|z-w|^{-d-\alpha} dw\nonumber\\
&\le c_2\|f\|_{\infty} z_d^{\beta_1}(|\log z_d|^{\beta_3}\vee 1)\int_{\R^d_+\setminus D(1,1)}\frac{1}{|w|^{d+\alpha+\beta_1}} 
\big(1+{\bf 1}_{|w|\ge1}(\log|w|)^{\beta_3}\big)
 dw\\
 &\le c_2 \|f\|_{\infty}z_d^{\beta_1} |\log z_d|^{\beta_3}\int_{\R^d_+\setminus D(1,1)}\frac{ \big(1+{\bf 1}_{|w|\ge1}(\log|w|)^{\beta_3}\big)}{|w|^{d+\alpha+\beta_1}} dw.
\end{align*}
Hence, by
 using the L\'evy system formula and \normal 
Lemma \ref{l:upper-bound-for-integral}, for $x\in D(2^{-5}, 2^{-5})$,
\begin{align}\label{e:bh-decay-3}
&\E_x\left[f(Y_{\tau_{U}}); Y_{\tau_{U}}\notin D(1,1)\right] \nn \\
& = \E_x \int_0^{\tau_U} \int_{\R^d_+\setminus D(1,1)} f(Y_t) 
 J(Y_t,w)\, dw \bk
\nn\\
&\le c_3 \|f\|_{\infty} x_d^{\alpha-1} \int_{\R^d_+\setminus D(1,1 )}\frac{ \big(1+{\bf 1}_{|w|\ge1}(\log|w|)^{\beta_3}\big)}{|w|^{d+\alpha+\beta_1}} dw\nn  \\
&\le c_4 \|f\|_{\infty} x_d^{\alpha-1}.
\end{align}
The claim of the lemma follows by using \eqref{e:bh.decay-2} and \eqref{e:bh-decay-3} in \eqref{e:bh-decay-1}. \qed

\begin{prop}\label{p:h-decay}
There exists $C>0$ such that for every $f:\R^d_+\to [0,\infty)$ which is regular harmonic 
in $U=D(1/2,1/2)$, it holds that
$$
f(x)\le C \left(\frac{f(w)}{w_d^{\alpha-1}}\right)   x_d^{\alpha-1}, \quad w, x\in D(2^{-3}, 2^{-3}).
$$
In particular, $f$ vanishes continuously $\{x=(\wt{x},0)\in \partial \R^d_+: |\wt{x}|< 2^{-3}\}$. 
\end{prop}
\pf For any $k\in \N$ define
$$ 
f_k(x):=\E_x[f(Y_{\tau_U})\wedge k]=\E_x[(f\wedge k)(Y_{\tau_U})]  .
$$
Then $f_k$ is regular harmonic in $U$ and bounded in $\R^d_+$. By Lemma \ref{l:bh-decay}, 
$f_k(x)\le c_4 \|f\wedge k\|_{\infty} x_d^{\alpha-1}$. 
Thus $f_k$ vanishes continuously  at the boundary. By the boundary Harnack principle, Theorem \ref{t:BHP},   there is $C>0$ such that for every $k\in \N$
$$
\frac{f_k(x)}{f_k(w)}\le C \frac{x_d^{\alpha-1}}{w_d^{\alpha-1}}, \qquad x, w\in D(2^{-3}, 2^{-3}).
$$
By letting $k\to \infty$ we get that
$$
\frac{f(x)}{f(w)}\le C \frac{x_d^{\alpha-1}}{w_d^{\alpha-1}}, \qquad x, w\in D(2^{-3}, 2^{-3}).
$$
\qed

\section{Estimates on Green functions and Potential}\label{s:EGP}

\subsection{Green function estimates}
 By following \cite[Section2]{KSV21} we see that there exists a symmetric function $ G(x,y)$, $x,y\in \R^d$ such that for all Borel functions $f:\R^d\to [0,\infty)$ we have that
$$
\E_x\int_0^{\zeta}f(Y_t)\, dt =\int_{\R^d_+}G(x,y)f(y)\, dy.
$$
Moreover, since $Y$ is transient, $G(x,y)$ is not identically infinite, and as in \cite[Proposition 2.2]{KSV21}, we can conclude that $G(\cdot, \cdot)$ is lower-semicontinuous in each variable and finite outside the diagonal. Further, for every $x\in \R^d_+$, $G(x, \cdot)$ is harmonic with respect to $Y$ in $\R^d_+\setminus\{x\}$ and regular harmonic in $\R^d_+\setminus B(x, \epsilon)$ for every $\epsilon >0$. The function $G$ enjoys the following scaling property (proved as in \cite[Proposition 2.4]{KSV21}): For all $x,y\in \R^d_+$, $x\neq y,$
\begin{equation}\label{e:green-scaling}
G(x,y)=G\left(\frac{x}{|x-y|}, \frac{y}{|x-y|}\right)|x-y|^{\alpha-d}\, .
\end{equation}

Choose a $p \in (\alpha-1, \alpha)$ and let
$
\kappa(x)
=C(\alpha, p,  \sB)x_d^{-\alpha}.
$
Let $Y^{\kappa}$ be the subprocess of $Y$ with killing potential $\kappa$ so that the corresponding Dirichlet form  is
$ \EE^{\kappa}(u,v)=\normal \EE(u,v)+\int_{\R^d_+} u(x)v(x)\kappa(x) dx$.
Let $G(x, y)$ and $G^{\kappa}(x, y)$  be the Green functions of $Y$ and $Y^{\kappa}$ respectively. Then $G(x, y)\ge G^{\kappa}(x, y)$. Now the
following result follows immediately from \cite[Proposition 4.1]{KSV21}.

\begin{prop}\label{p:green-lower-bound}
For any $C_1>0$, there exists a constant $C_2>0$ 
such that for all $x,y\in \R^d_+$ satisfying 
$|x-y|\le C_1(x_d\wedge y_d)$, it holds that
$$
G(x,y)\ge C_2|x-y|^{-d+\alpha}.
$$
\end{prop}

From now on we assume $d > (\alpha+\beta_1+\beta_2)\wedge 2$. 
In \cite[Section 4.2]{KSV21}, the killing function plays no role. Repeating the argument  leading to \cite[Proposition 4.6]{KSV},  we get

\begin{prop}\label{p:green-upper-bound}
There exists a constant $C>0$ 
such that for all $x,y\in \R^d_+$ satisfying 
$|x-y|\le x_d\wedge y_d$, it holds that
$$
G(x,y) \le C |x-y|^{-d+\alpha}.
$$
\end{prop}

Using Proposition \ref{p:h-decay}, 
we can combine Proposition \ref{p:green-upper-bound}
with Theorem \ref{t:carleson} to get the following result, which is key for us to get sharp
two-sided Green functions estimates.

\begin{prop}\label{p:gfcnub}
There exists a constant $C>0$ 
such that for all $x, y\in \R^d_+$,
\begin{align}
\label{e:upper}
G(x,y) \le C |x-y|^{-d+\alpha}.
\end{align}
\end{prop}

\pf
It follows from Proposition \ref{p:green-upper-bound} that there exists $c_1>0$ such that
$G(x, y)\le c_1$ for all $x, y\in \R^d_+$ with $|x-y|=1$ and $x_d\wedge y_d>1$.
By Theorem \ref{t:uhp}, for any $c_2>0$, there exists $c_3>0$ such that 
$G(x, y)\le c_3$ for all $x, y\in \R^d_+$ with $|x-y|=1$ and $x_d\wedge y_d>c_2$.
Now by Theorem \ref{t:carleson}, we see that there exists $c_4>0$ such that
$G(x,y)  \le c_4$ for all $x, y\in \R^d_+$ with $|x-y|=1$.
Therefore, by  \eqref{e:green-scaling}, we have
$$
G(x,y) \le C |x-y|^{-d+\alpha}, \quad x, y\in \R^d_+. 
$$
\qed

Now we prove the two-sided Green function estimates.

\noindent
{\bf Proof of Theorem \ref{t:Green}.}
The scaling property and the invariance property of the half space under scaling imply that in order to prove Theorem \ref{t:Green}  we only need to show that for all $x,y\in \R^d_+$ satisfying $|x-y|=1$, 
\begin{equation}\label{e:Green2}
C^{-1} \left(x_d \wedge 1 \right)^{\alpha-1}\left({y_d} \wedge 1 \right)^{\alpha-1} \le G (x,y) \le C  \left(x_d \wedge 1 \right)^{\alpha-1}\left({y_d} \wedge 1 \right)^{\alpha-1}.
\end{equation}
By scaling, Theorem \ref{t:uhp},  and Propositions \ref{p:green-lower-bound} and
\ref{p:gfcnub},
we only need to show \eqref{e:Green2} for $x_d \wedge y_d \le 2^{-3}$ and $|x-y|=1$.

We now assume that  $|x-y|=1$ and let  $x_0=(\wt x, 2^{-3})$ and $y_0=(\wt y, 2^{-3})$. We first note that Proposition \ref{p:h-decay}, together with scaling, clearly implies that for all $x\in \R^d_+$, $y\mapsto G(x,y)$ vanishes at the boundary of $\R^d_+$. 

Suppose that $y_d \ge 2^{-3}$. Then, 
by Theorem \ref{t:uhp}, and Propositions \ref{p:green-lower-bound} and \ref{p:gfcnub}, 
we have $G(x_0,y) \asymp c >0$.
Thus by Theorem \ref{t:BHP}, 
\begin{align}
\label{e:Gs1}
G(x,y) \asymp  G(x_0,y) (x_d/2^{-3})^{\alpha-1}  \asymp x_d^{\alpha-1}.
\end{align}

Suppose that $y_d < 2^{-3}$. Then 
by Theorem \ref{t:BHP} and \eqref{e:Gs1}
\begin{align}
\label{e:Gs2}
G(x,y) \asymp  G(x,y_0) (y_d/2^{-3})^{\alpha-1}  \asymp x_d^{\alpha-1} y_d^{\alpha-1}.
\end{align}
\eqref{e:Gs1} and \eqref{e:Gs2} imply that 
\eqref{e:Green2} hold for $x_d \wedge y_d \le 2^{-3}$ and $|x-y|=1$.
\qed

\subsection{Estimates on  Potentials}\label{ss:lb}

Recall that we have assumed $d > (\alpha+\beta_1+\beta_2)\wedge 2$. 
Let $G^{B(w,R)\cap \R^d_+}(x,y)$ be 
the Green function of the process $Y$ killed upon exiting 
$B(w,R)\cap \R^d_+$, $w\in \partial \R^d_+$.

For any $a>0$, let $B^+_a:=B(0, a)\cap \R^d_+$. 
Using Theorem \ref{t:Green} and the formula
$$
G^{B^+_1}(y, z)=G(y, z)-\E_y[G(Y_{\tau_{B^+_1}}, z)], 
$$
the proof of next result is standard. For example, see \cite[Lemma 5.1]{KSV21}
\begin{lemma}\label{l:GB_1}
For any $\eps \in (0, 1)$ and $M>1$, there exists a constant 
$C>0$ such that for all 
$y,z \in B^+_{1-\eps}$ with $|y-z|  \le M(y_d\wedge z_d)$.
$$
G^{B^+_1}(y, z)\ge C |y-z|^{-d+\alpha}.
$$
\end{lemma}

By Theorems \ref{t:BHP} and \ref{t:uhp}  we have for any $r>0$ and $x\in \R^d_+$ with $x_d<r/2$,
$$
\P_x(Y_{\tau_{D_{\wt x}(r,r)}} \in D_{\wt x}(r, 4r) \setminus D_{\wt x}(r, 3r))
\ge c x^{\alpha-1}.
$$
Using this and Lemma \ref{l:GB_1}, the proofs of next two lemmas are identical to 
those of \cite[Lemmas 5.2 and 5.3]{KSV21}.
\begin{lemma}\label{p:GB_1}
For every $\eps \in (0, 1/4)$ and $M, N>1$, there exists a constant 
$C_>0$ such that for 
all $x,z \in B^+_{1-\eps}$ with $x_d \le z_d$ 
satisfying $x_d/N \le |x-z|\le M z_d $, it holds that
$$
G^{B^+_1}(x, z)\ge Cx^{\alpha-1}_d|x-z|^{-d+1}.
$$
\end{lemma}

\begin{lemma}\label{l:GB_4}
For every $\eps \in (0, 1/4)$ and $M \ge 40/\eps$, there exists a constant 
$C>0$ such that 
for all $x,z \in B^+_{1-\eps}$ with $x_d \le z_d$ 
satisfying $|x-z|\ge M z_d $, it holds that
$$
G^{B^+_1}(x, z)\ge C
x^{\alpha-1}_dz^{\alpha-1}_d|x-z|^{-d-\alpha+2}.
$$
\end{lemma}
Combining the above results with scaling, we get 
\begin{thm}\label{t:GB}
For any $\eps \in (0, 1/4)$, 
there exists a constant $C>0$ such that for all 
$w \in \partial \R^d_+$,  $R>0$ and $x,y \in 
B(w, (1-\eps)R)\cap \R^d_+$, it holds that
$$
G^{B(w, R)\cap \R^d_+}(x, y)\ge C
\left(\frac{x_d}{|x-y|}  \wedge 1 \right)^{\alpha-1}\left(\frac{y_d}{|x-y|}  \wedge 1 \right)^{\alpha-1} \frac{1}{|x-y|^{d-\alpha}}.
$$
\end{thm}

\begin{prop}\label{p:bound-for-integral-new} 
For any 
$\wt{w} \in \R^{d-1}$, any Borel set $D$ 
satisfying $D_{\wt{w}}(R/2,R/2) \subset D  \subset D_{\wt{w}}(R,R)$ and any 
$x=(\wt{w}, x_d)$ with $x_d \le R/10$,
\begin{equation}\label{e::bound-for-integral-new} 
\E_x \int_0^{\tau_D}(Y_t^d)^{\gamma \normal}\, dt  = \int_D G^D(x,y)y_d^{\gamma \normal}\, dy  \asymp  
\begin{cases} 
R^{ \gamma  +1}x_d^{\alpha-1}, & \gamma \normal>-1,\\ 
x_d^{\alpha-1}\log(R/x_d), &   \gamma =-1, \\ 
x_d^{\alpha+ \gamma \normal}, &-{\alpha}<  \gamma <-1, 
\end{cases}
\end{equation}
where the comparison constant is independent of $\wt{w} \in \R^{d-1}$, $D$, $R$ and $x$.
\end{prop}
\pf
Let  $D$ be a Borel set satisfying $D(R/2,R/2) \subset D  \subset D(R,R)$. 
By Theorems \ref{t:Green} and \ref{t:GB}, we have that for all $x \in D$ 
\begin{align}
&\int_D G^D(x,y)y_d^{\gamma \normal}\, dy \le \int_D G(x,y)y_d^{\gamma \normal}\, dy\nn\\
& \le
c_1 \int_{D(R, R)}  \left(\frac{x_d}{|x-y|}  \wedge 1 \right)^{\alpha-1}\left(\frac{y_d}{|x-y|}  \wedge 1 \right)^{\alpha-1} \frac{ y_d^{\gamma}dy}{|x-y|^{d-\alpha}} ,  
\label{e:e:pu1}
\end{align}
and for  $x=(\wt{0}, x_d)$ with $x_d \le R/10$
\begin{align}
& \int_D G^D(x,y)y_d^{\gamma \normal}\, dy 
 \ge \int_{B^+_{R/2}}  y_d^{\gamma} 
G^{B^+_{R/2}}(x, y)\, dy 
\ge c_2\int_{D(R/5, R/5)}   y_d^{\gamma} 
G^{B^+_{R/2}}(x, y)\, dy\nn\\
&\ge \int_{D(R/5, R/5)}  \left(\frac{x_d}{|x-y|}  \wedge 1 \right)^{\alpha-1}\left(\frac{y_d}{|x-y|}  \wedge 1 \right)^{\alpha-1} \frac{ y_d^{\gamma}dy}{|x-y|^{d-\alpha}}.
\label{e:e:pd1}
\end{align}
We now apply \cite[Lemma 3.3 and Theorem 3.4]{AGV} and their proofs to \eqref{e:e:pu1} and \eqref{e:e:pd1} and get \eqref{e::bound-for-integral-new} 
\qed

\begin{corollary}
For all $x \in \R^d_+$,
$$
\E_x\int_0^{\zeta} (Y_t^d)^{ \gamma \normal}\, dt =   \int_{\R^d_+}G(x,y)y_d^{\gamma \normal} dy
 \asymp  
\begin{cases} 
\infty &  \gamma \normal \ge -1  \text{ or }  \gamma \normal \le -{\alpha},\\ 
x_d^{\alpha+ \gamma \normal}, & -\alpha<   \gamma \normal<-1.
\end{cases}
$$
In particular,  for all $x \in \R^d_+$,
$\E_x[\zeta]= \infty$.
\end{corollary}
\pf
When
$ \gamma>-\alpha $,
the result follows by letting $R \to \infty$ in Proposition \ref{p:bound-for-integral-new}.
If $ \gamma \le -\alpha $,
by Theorems \ref{t:Green}, for all $x \in \R^d_+$
 \begin{align*}
&  \int_{\R^d_+}G(x,y)y_d^{\gamma \normal} dy \ge c_1 \int_{D_{\wt{x}}(x_d, x_d/2)} \left( \frac{x_d}{|x-y|}  \right)^{\alpha-1}\left(\frac{y_d}{|x-y|} \right)^{\alpha-1} \frac{ y_d^{\gamma}dy}{|x-y|^{d-\alpha}}\\
&= c_1x_d^{\alpha-1} \int_{D_{\wt{x}}(x_d, x_d/2)} 
\frac{ y_d^{\gamma+\alpha-1}dy}{|x-y|^{d+\alpha-2}} \ge c_2x_d^{1-d}  
\int_{D_{\wt{x}}(x_d, x_d/2)} y_d^{\gamma+\alpha-1}dy=\infty.
\end{align*}
\qed

\vskip 0.1truein

\parindent=0em

{\bf Panki Kim}

Department of Mathematical Sciences and Research Institute of Mathematics,

Seoul National University, Seoul 08826, Republic of Korea

E-mail: \texttt{pkim@snu.ac.kr}

\bigskip

{\bf Renming Song}

Department of Mathematics, University of Illinois, Urbana, IL 61801,
USA

E-mail: \texttt{rsong@math.uiuc.edu}

\bigskip

{\bf Zoran Vondra\v{c}ek}

Department of Mathematics, Faculty of Science, University of Zagreb, Zagreb, Croatia,

Email: \texttt{vondra@math.hr}

\end{document}